\documentclass[11pt,leqno]{amsart}

\usepackage[sc]{mathpazo}
\usepackage[T1]{fontenc}

\usepackage[a4paper,hmargin=2.8cm,vmargin=4cm]{geometry}
\usepackage{amsfonts,amssymb,amscd,amstext}
\usepackage{mathtools}
\mathtoolsset{showonlyrefs=true}
\usepackage[utf8]{inputenc}
\usepackage[colorlinks=true,linkcolor=blue,citecolor=red]{hyperref}
\usepackage{braket,relsize,nccmath}
\newtheorem{proposition}{Proposition}
\newtheorem{theorem}[proposition]{Theorem}
\newtheorem*{th1}{Theorem A}
\newtheorem*{th2}{Theorem B}
\newtheorem{corollary}[proposition]{Corollary}
\newtheorem{remark}{Remark}

\newfont{\bb}{msbm10 at 12pt} 

\def\pf{{\textit {Proof :}}}

\def\C{{\mathbb C}}

\def\S{\hbox{\bb S}}

\def\ric{{\rm Ric}}
\def\div{{\rm div}}
\def\calp{{\mathfrak p}}

\newcommand{\bal}{\begin{align}}      \newcommand{\eal}{\end{align}}
\newcommand{\ba}{\begin{array}}      \newcommand{\ea}{\end{array}}
\newcommand{\bc}{\begin{center}}     \newcommand{\ec}{\end{center}}
\newcommand{\be}{\begin{enumerate}}  \newcommand{\ee}{\end{enumerate}}
\newcommand{\beq}{\begin{eqnarray}}  \newcommand{\eeq}{\end{eqnarray}}
\newcommand{\beQ}{\begin{eqnarray*}} \newcommand{\eeQ}{\end{eqnarray*}}
\newcommand{\bi}{\begin{itemize}}    \newcommand{\ei}{\end{itemize}}
\newcommand{\bt}{\begin{tabular}}    \newcommand{\et}{\end{tabular}}
\newcommand{\bdm}{\begin{displaymath}} \newcommand{\edm}{\end{displaymath}}

\let\pa=\partial

\newcommand{\lto}{\ensuremath{\longrightarrow}}

\def\qed{\hfill{q.e.d.}\smallskip\smallskip}

\begin{document}

\title[A compact approach to the positive mass of Brown-York mass]{A compact approach to the positivity of Brown-York's mass and its relation with the Min-Oo conjecture, Yau's Problem \#100 and  rigidity of hypersurfaces}

\author{Sebasti{\'a}n Montiel}
\address[Montiel]{Departamento de Geometr{\'\i}a y Topolog{\'\i}a\\
Universidad de Granada\\
18071 Granada \\
Spain}
\email{smontiel@ugr.es}

\begin{abstract}
In this paper, we fix some mistakes in previous versions and we elaborate a simple spinorial proof for the positivity of the Brown-York mass avoiding the glueing of infinite pieces along the  boundary of the compact manifold. Also, we will not need either to solve difficult PDE's in a $C^1$ context or the previously known positivity of the ADM mass. Moreover, we will prove versions of the Min-Oo conjecture, of famous Yau's Problem \#100 about embedded minimal surfaces in the unit three-sphere and of the Fenchel inequality.
\end{abstract}

\thanks{Research partially supported by a Junta de Andaluc\'ia grant FQM-325}



\date{September, 2024}        

\maketitle \pagenumbering{arabic}

\begin{flushright}
{\em{\tiny Aní lah ve-att li..., `od?}}
\end{flushright}
\vspace{.5cm}

\section{Introduction}
Shi and Tam  proved in \cite[Theorem 4.1]{ST1} the expected positivity for the Brown-York (BY) \cite{BY} semi-local mass. They based, among other things, on a $C^1$version of the previously proved positivity of the  Arnowitt, Deser, Misner  (ADM) mass \cite{ADM}. These tools are perhaps the most used and useful ones in General Relativity to determine the energy enclosed in a compact domain $\Omega$ of a three-dimensional spacelike slice of a given spacetime. If we forget the physical meaning of this result, it has been laid bare in plenty works by several authors (see \cite{GHHP,LiuY1,LiuY2,MiST,Mu,ST2,WaY1,WaY2}) that we can interpret the positivity of the BY mass in a  purely Riemannian context. In fact, if $\Omega$ is a three-dimensional compact and connected Riemannian  manifold with non-negative scalar curvature and 
strictly convex boundary $\Sigma$, then 
$$
\int_{\Sigma}H_0 -\int_{\Sigma}H\ge 0,
$$
where $H$ is the inner mean curvature of $\Sigma$
in $\Omega $ and $H_0$ is the Euclidean inner mean
curvature of $\Sigma$ in ${\mathbb R}^{3}$ given by the (unique up to rigid motions) Weyl convex embedding of $\Sigma$ into ${\mathbb R}^3$. Pogorelov \cite{Po} and Nirenberg \cite{Ni} showed  independently the existence and uniqueness of this embedding for abstract convex surfaces.  Moreover, the equality holds for some boundary connected component $\Sigma$ if and only if $\Omega$ is a Euclidean domain. The inventive proof by Shi and Tam uses two fundamental facts. First, following an idea by Bartnik \cite{Bar}, the construction of a suitable infinite asymptotically flat extension of $\Omega$ with non-negative scalar curvature by glueing one infinite piece along $\Sigma$ in ${\mathbb R}^3$ in order it has the same mean curvature as $\Sigma$ has in $\Omega$. This construction is equivalent to solve a non-linear parabolic equation in a $C^1$ context, because this is the degree of differentiability after the glueing. The second point is to see that the difference of integrals whose limit defines the Brown-York mass, given, as we have already seen above by$$
\int_{\Sigma_r}H_0-\int_{\Sigma_r}H_r,$$ 
where $\Sigma_r$ is an expansion to infinity of the original boundary $\Sigma$ into ${\mathbb R}^3$, $H_0$ is the mean curvature with respect to the Euclidean metric and $H_r$ is the mean curvature with respect to the Bartnik metric, goes in a non decreasing way to the ADM mass of the  asymptotically flat manifold built before, and, so, one can get the non-negativity of this mass. This approach was successfully used to prove the positivity of other quasi local masses proposed by Liu and Yau \cite{LiuY1,LiuY2} and Wang and Yau \cite{WaY1}. 

Along these comments one can see that positivity of quasi local masses and rigidity of compact manifolds with non-empty boundary are different, although closely related, aspects of almost identical questions. 

On one hand, we know that Witten, by using a spinorial structure on $\Omega$ (the topological obstructions to support a  spin structure are always satisfied on  oriented manifolds with dimension three), elaborated a proof for the theorem of the positivity of the ADM mass defined on complete asymptotically flat with positive scalar curvature \cite{Wi,He1,He2} which remarkably simplifies the original one  by Schoen and Yau. The fundamental difference between both proofs is that the original one uses the existence and behaviour of minimal surfaces in $\Omega$ and the second one makes use of certain spinors defined on $\Omega$ and compare their conduct with that of the restriction to the compact surfaces embedded in the bulk manifold. 

In another order of things, following an old clever observation by Jean-Pierre Bourguignon, the author, his colleagues B\" ar \cite{Ba2}, Hijazi \cite{Hi}, Raulot \cite{HMRa2} {\em et al.}\ 
  began to introduce the spin machinery and the Dirac operator in the study of hypersurfaces, with the aim of simplify the proofs of some well-known theorems and, maybe, obtain other new ones. Before that, several previous attempts in this direction had been tried, fundamentally by Trautman \cite{Tr} and Bure{\v s} \cite{Bur}. In the last years, a lot of works about submanifolds theory make use of the comparison between the Dirac operator $D$ acting on spinor fields on $\Omega$ and the Dirac  operator ${\bf D}$ of the  induced  spin structrure on the boundary $\Sigma= \pa\Omega$ which acts on the corresponding restrictions.  This change of the Laplacian operator for the Dirac one has its advantages and its drawbacks. The first ones are focused on the fact that the difference between the extrinsic  and intrinsic Dirac operators $D$ and ${\bf D}$ acting on the spinor fields of $\Sigma$ is controlled by the scalar curvature of $\Omega$ and the mean curvature of $\Sigma$. Instead, it is well-known that, in the case of the Laplacian, this same difference is governed by the Ricci curvature of $\Omega$ and the whole of the second fundamental form of $\Sigma$ in  $\Omega$.  As for the downsides, first, we have that not all the orientable Riemannian manifolds can support a spin structure. But the corresponding topological obstruction is the vanishing of its second Whitney  class $w_2(\Omega)$. But, in fact, this is quite weak condition. It is accomplished by all the usual ambient spaces where geometers of submanifolds work.   Let us add to that, that moreover,  $w_2(\Omega)=0$ for all the (orientable) three-dimensional Riemannian manifolds.   Secondly, the restriction of a function defined on any manifold to any of its submanifolds, without restrictions on its codimension, is obviously another function. However, the restriction of the spinor bundle defined on a given spinorial Riemannian manifold is not, in general, the spinor bundle induced on its submanifolds. The reader can see in \cite{Ba2} that this restriction involves the structure of the normal bundle of the submanifold, or, to give a stranger example, it has been shown in \cite{HMU} that the spinorial fields on a K\"ahler manifolds may be identifiable with certain differentiable exterior forms on its Lagrangian submanifolds.   

Whe will remember as, since sometime ago, this spinorial approach, applied to the theory of submanifolds,  allowed us to estimate the eigenvalues of the Dirac operator ${\bf D}$ of a compact surface $\Sigma$ enclosing a domain $\Omega$ in terms, exclusively, of the mean curvature of the boundary (see \cite{HMZ1,HMZ2}). Moreover, from these estimates,  we  got an almost trivial proof for the Alexandrov theorem for the embedded hypersurfaces with constant mean curvature \cite{HMZ2}, we consider and solved the Alexandrov question for certain spacelike slices of a the Lorentz-Minkowski spacetime (\cite{HMRa1}), we related the first non-null eigenvalue $\lambda_1({\bf D})$  of an embedded compact hypersurface of the Euclidean space with the first eigenvalue $\nu_1$ of the corresponding Steklov problem on $(\Omega,\Sigma)$  and we gave an answer to the a {\em flat version} of the Min-Oo's conjecture \cite{HMZ2}.   In fact, in the flat case, our  answer to the corresponding {\em flat} Min-Oo conjecture was also independently obtained by Miao \cite{Mi1} (see Remark 1). In our case, this answer is not but a mere consequence of our estimate for $\lambda_1({\bf D})$.

In this article, we will fixed some grave mistakes in the proofs in our recent paper \cite{Mon2}. We discovered these mistakes thanks to a question posed to us by one  of the referees of the journal where the paper was submitted and it was also pointed out to us by our colleague and friend Christian B\"ar.  Once we have done this, the paper is dedicated to expose a compact approach to the aforementioned Shi-Tam result about the positivity of the Brown-York mass. Besides to provide an easy spinorial proof, we generalised it to mean-convex bodies. In the aforementioned paper, we emphasised the physical side of things and, in some remarks, we scarcely indicated in passing the strong relation between them and some results on the embedded (hyper)surfaces, especially in Euclidean spaces and spheres. Now, we will fix \cite{Mon2} with geometer eyes and show how the spinor estimates found there, with their corresponding corrections, provide simplified proofs, improve famous theorems and give us  new nice results. 

So, the aim of writing this paper is twofold. On one hand, we want to bring together the main new results and some improvements of old ones that we had already obtained in the study of hypersurfaces. We will see how apparently unrelated theorems can be deduced from a same source. On the other hand, we will obtain, but using again from the same spinorial mechanism, new results about (hyper)surfaces in Euclidean and spherical ambient spaces. For example, a spinorial version of the famous \#100 problem posed by Yau for the Laplacian in \cite{Y2}. 

\vspace{0.5cm}
 \noindent
{\em Acknowledgments:} The author  would like to  thank Prof.\  J.-P. Bourguignon for his far-sighted encouragement so that the geometers include the spinorial geometry among their mathematical tools with the same intensity as we made use of the Laplacian, and so that we should avoid our resistance to the alleged difficulties inherent to the first order differential operators.  Also, our colleagues Oussama Hijazi (Universit\'e de Lorraine) and Simon Raulot (Universit\'e de Rouen), who have worked with us, since the beginning of century, following  the valuable piece of advice of Bourguignon (though the pioneers were Trautman, Anghel and B\^ures  \cite{Tr,Ang,Bur} {\em et al.}). Equally, we thank Christian B\"ar (Potsdam Universit\"at) for valuable e-conversations and his papers' clarity. In another order of things, we want to thank as well the Judge of the Criminal Court no.\! 3 of Granada. He was empathetic enough to grant us the necessary  time and the  means to work in and complete the redaction of this article.

\section{Riemannian spin manifolds and hypersurfaces}\label{spin} 
At the risk of being too long and taking in mind that the spinorial geometry is little used in theory of submanifolds, we will start by briefly recalling the fundamental properties of spinorial structures  which will utilize, and how to manage them to obtain nice results in hypersurfaces theory. 

Let  $\Omega$ be an $(n+1)$-dimensional spinorial (oriented) Riemannian manifold  and let $\Sigma$ be an orientable submanifold with dimension $n$ immersed as a hypersurface 
in $\Omega$. We  denote by $\langle\,,\rangle$ the scalar product on $\Omega$, by $\nabla$ its corresponding Levi-Civita connection and by $\gamma:\C\ell(\Omega)\lto {\rm End}_{\C}(\S\Omega)$ its Clifford multiplication, where  $\S\Omega$ is the corresponding spinorial bundle of the ambient space,  a complex vector bundle of rank $2^
{\left[\frac{n+1}{2}\right]}$ whose fibers are preserved by the irreducible representation of the Clifford algebras constructed over the tangent spaces of $\Omega$. When the dimension $n+1$ is even, we have the orthogonal chirality standard decomposition
\begin{equation}\label{chi-de}
\S\Omega=\S\Omega^+\oplus\S\Omega^-,
\end{equation}
where the two direct summands are respectively the $\pm 1$-eigenspaces of the endomorphism 
$\gamma(\omega_{n+1})$, with $\omega_{n+1}=i^{\left[\frac{n+2}{2}\right]}e_1\cdots e_{n+1}$ is the so-called  complex volume form of $\Omega$. It is well-known (see \cite{LM}) that there are, on the rank $2^{\left[\frac{n+1}{2}\right]}$  complex spinor bundle $\S\Omega$, a natural Hermitian metric denoted also by $\left<,\right>$  and a spinorial Levi-Civita connection  represented by $\nabla$ as well, all of them compatible with both the original metric  on the basis  $\Omega$ and with the Clifford multiplication $\gamma$. This compatibility means
\beq
 &X\left<\psi,\phi\right>=\left<\nabla_X\psi,\phi\right>+\left<\psi,\nabla_X\phi\right>& \label{nabla-metric} \\
&\nabla_X\left(\gamma(Y)\psi\right)=\gamma(\nabla_XY)\psi+\gamma(Y)\nabla_X\psi& \label{nabla-gamma}
\eeq
for any tangent vector fields $X,Y\in\Gamma(T\Omega)$ and any spinor fields $\psi,\phi\in\Gamma
(\S\Omega)$. Moreover, with respect to this Hermitian product on $\S\Omega$, Clifford multiplication by vector fields is skew-Hermitian or equivalently
\begin{equation}\label{skew}
\left<\gamma(X)\psi,\gamma(X)\phi\right>=|X|^2\left<\psi,\phi\right>.
\end{equation}
Since the complex volume form $\omega_{n+1}$ is parallel with respect to the spinorial Levi-Civita connection, when $n+1=\dim \Omega$ is even, the chirality decomposition 	\eqref{chi-de} is preserved by $\nabla$. Moreover, from \eqref{skew}, one sees that this decomposition is also orthogonal.

Let us consider the Dirac first order differential operator $D$ on the spinorial manifold $\Omega$, acting on the sections of ${\mathbb S}\Omega$,
$$
D=\sum_{i=1}^{n+1}\gamma(e_i)\nabla_{e_i},$$
where $\{e_1,\dots,e_{n+1}\}$ is any local orthonormal frame in $T\Omega$. When $n+1=\dim \Omega$ is even, $D$ interchanges the chirality subbundles $\S\Omega^\pm$. Its relation with the much more known and used rough operator $\Delta$, of second order, is given by 
\begin{equation}\label{Dq}
{ D}^2=-\Delta+\frac{R}{4}=\nabla^*\nabla+\frac{R}{4},
\end{equation}
where $R$ is the scalar curvature of $\Omega$.

On the orientable hypersurface $\Sigma$, we choose the orientation which makes an orthonormal basis $\{e_1,\dots,e_n\}$ positive just when the completed orthonormal basis $\{e_1,\dots,e_n,N\}$ is positive on $\Omega$, where $N$ is the unit normal vector at each point of $\Sigma$ associated to the orientation of the ambient space. We consider the induced Riemannian metric from $\Omega$, also denoted by $\left<\,,\right>$ and its Levi-Civita connection $\nabla^{\Sigma}$. It is a standard fact that both connections are related by means of  the so-called Gauss and Weingarten equations, where appear  the shape operator $A$ of the hypersurface $\Sigma$ corresponding to the unit normal field $N$. As the normal bundle of the  hypersurface is trivial, the Riemannian manifold $\Sigma$ has also a {\em restricted} spin manifold and so we will have the corresponding rank $2^{\left[\frac{n}{2}\right]}$ spinor bundle $\S\Sigma$, the Clifford multiplication $\gamma^{\Sigma}$, the spinorial Levi-Civita connection $\nabla^{\Sigma}$ and the intrinsic Dirac operator $D^{\Sigma}$. It is not difficult to show (see \cite{Ba2,BFGK,Bur,Tr,Mon1}) that the restricted Hermitian bundle$$
{\bf S}=\S\Omega_{|\Sigma}$$
can be identified with the intrinsic Hermitian spinor bundle $\S\Sigma$, provided that $n=\dim \Omega$ is odd. Instead, if $n=\dim \Omega$ is even, the restricted bundle
${\bf S}$ can be identified with the sum $\S\Sigma^+\oplus\S\Sigma^-$. For any spinor field $\psi\in\Gamma({\bf S})$ on the boundary hypersurface $\Sigma$ and any vector field $X\in\Gamma(T\Sigma)$, we define on the restricted bundle ${\bf S}$, a Clifford multiplication $\gamma^{{\bf S}}$ and a connection $\nabla^{{\bf S}}$ by
\beq\label{nablaS}
& \gamma^{{\bf S}} (X)\psi=\gamma(X)\gamma(N)\psi & \\
& \nabla^{{\bf S}}_X\psi=\nabla_X\psi-\frac{1}{2}\gamma^{{\bf S}}(AX)\psi=\nabla_X\psi-
\frac{1}{2}\gamma(AX)\gamma(N)\psi\,.\nonumber
\eeq
 Then it easy to see that $\gamma^{{\bf S}}$ and $\nabla^{{\bf S}}$ correspond respectively to $\gamma^{\Sigma}$ and $\nabla^{\Sigma}$, for $n$ odd, and to $\gamma^{\Sigma}\oplus -\gamma^{\Sigma}$ and $\nabla^{\Sigma}\oplus\nabla^{\Sigma}$, for $n$ even. Then, $\gamma^{{\bf S}}$ and $\nabla^{{\bf S}}$ satisfy the same compatibilty relations  
\eqref{nabla-metric}, \eqref{nabla-gamma} and \eqref{skew}, together
with the following additional identity 
\begin{equation}\label{Supercom}
\nabla^{{\bf S}}_X\left(\gamma(N)\psi\right)=\gamma(N)\nabla^{{\bf S}}_X\psi.
\end{equation}
As a consequence, the hypersurface Dirac operator ${\bf D}$ acts on smooth sections $\psi\in\Gamma({\bf S})$ as follows
\begin{equation}\label{bD}
{\bf D}\psi =\sum_{j=1}^{n}\gamma^{{\bf S}}(u_j)\nabla^{{\bf S}}_{u_j}\psi=
\frac{n}{2}H\psi-\gamma(N)\sum_{j=1}^{n}\gamma(u_j)\nabla_{u_j}\psi,
\end{equation}
where $\{u_1,\dots,u_{n}\}$ is a positive local orthonormal frame tangent to the boundary $\Sigma$ and $H=(1/n){\rm trace}\, A$ is its mean curvature function associated to $N$. So,
${\bf D}$ coincides with the intrinsic Dirac operator $D^{\Sigma}$ on the boundary, for $n+1$  odd, and with the pair $D^{\Sigma}\oplus -D^{\Sigma}$, for $n+1$ even. In the particular case where the field $\psi\in\Gamma({\bf S})$ is the restriction of a spinor field $\psi\in\Gamma(\Sigma)$ on $\Omega$, this means that
\begin{equation}\label{twoD}
{\bf D}\psi=\frac{n}{2}H\psi-\gamma(N)D\psi-\nabla_N\psi.
\end{equation}
Note that we always have the anticommutativity property
\begin{equation}\label{supercom}
{\bf D}\gamma(N)=-\gamma(N){\bf D}.
\end{equation}
So, when $\Sigma$ is compact, the spectrum of ${\bf D}$ is symmetric with respect to zero and coincides with the spectrum of $D^{\Sigma}$, for $n+1$ odd, and with ${\rm Spec}(D^{\Sigma})\cup -{\rm Spec}(D^{\Sigma})$, for $n+1$ even (see, fro example, \cite{HMRo}). Another well-known fact is that the spectrum Spec($\bf D$)  is a ${\mathbb Z}$-symmetric sequence$$
 -\infty\swarrow\cdots\le -\lambda_{k}\le\cdots\le-\lambda_{1}<\lambda_0=0<\lambda_1\le\cdots\le\cdots\le\lambda_k\le\cdots\nearrow
 +\infty,
  $$
where each eigenvalue is repeated according to its corresponding multiplicity and where $\lambda_0=0$ can have  multiplicity zero. This is because $D$ is an elliptic operator of order one which is self-adjoint due to the compacity of $\Sigma$ and  the fact that $\gamma(N)$ maps the eigenspace associated to $\lambda_k$ into that of $-\lambda_k$, for all
$k\in{\mathbb Z}$. 

\begin{remark}\label{highD}{\rm 
For higher codimensions $m\ge 1$, we recommend to read \cite{Ba2} for getting a correct understanding of the case. The restricted spin bundle ${\bf S}={\mathbb S}\Omega_{|\Sigma}$  is identifiable with ${\mathbb S}\Sigma\otimes  {\mathbb S}
\Sigma^\perp$, where ${\mathbb S} \Sigma^\perp$  is the rank $2^{[\frac{m}{2}]}$ spinor bundle built from the normal bundle of $\Sigma$ into $\Omega$, 
  unless $n$ and $m$ are both odd, in which case 
${\bf S}=\big({\mathbb S}\Omega_{|\Sigma}\otimes  {\mathbb S}\Omega^\perp_{|\Sigma}\big)\oplus \big(
{\mathbb S}\Omega_{|\Sigma}\otimes  {\mathbb S}\Omega^\perp_{|\Sigma}\big)$. Let $\nabla^{{\mathbb S}\Sigma}$ and  $\nabla ^{{\mathbb S}{\Sigma}^\perp}$  be the Levi-Civita connections on ${\mathbb S}\Sigma$ and ${\mathbb S}{\Sigma^
\perp}$.  We consider the connection $\nabla^{{\mathbb S}\Sigma}\otimes  \nabla^{{\mathbb S}\Sigma^\perp}$ given by ${\nabla ^{{\mathbb S}\Sigma}}\otimes Id+Id\otimes {\nabla^{{\mathbb S}\Sigma^\perp}}$. We mean it the product connection on ${\mathbb S}\Sigma\otimes {\mathbb S}{\Sigma^\perp}$ and also on $({\mathbb S}\Sigma\otimes {\mathbb S}{\Sigma^\perp})\oplus ({\mathbb S}\Sigma\otimes {\mathbb S}{\Sigma^\perp})$  if $n$ and $m$ are both odd. 
One conclude that
\begin{eqnarray}\label{cod}
&\gamma^{\bf S}=\gamma^{{\mathbb S}\Sigma}\otimes Id,&
\\&\nabla^{\bf S}= \nabla^{{\mathbb  S}\Sigma}   \otimes Id+Id\otimes 
\nabla^{{\mathbb S}\Sigma^\perp}  +\frac1{2}\sum^{n,m}_{i,j=1}\gamma(e_i)\gamma(\sigma(e_i, \cdot),\nonumber 
\end{eqnarray}
where $\sigma$ is the second fundamental form of the immersion of $\Sigma$ into $\Omega$. Furthermore, if ${\bf D}$ is the Dirac operator of the rank $2^{[\frac{n}{2}]+[\frac{m}{2}]}$ Clifford bundle ${\mathbb S}\Omega_{|\Sigma}$, in both cases, where $n$ and $m$ either are even or odd numbers, the Clifford bundles of type ${\mathbb S}\Sigma\otimes {{\mathbb S}\Sigma^\perp}$ are referred to as {\em the spinor bundle ${\mathbb S}\Sigma$ {\em twisted} by the vectorial bundle ${{\mathbb S}\Sigma^\perp}$}. In general, ${\mathbb S}\Omega_{|\Sigma}$, $m\ge 2$, is not a spinor bundle on $\Sigma$, but only a Clifford bundle. The rank of that type of bundles on an $n$-dimensional Riemannian manifold is greater than or equal to $2^{[\frac{n}{2}]}$ and the equality is attained
just by the spin bundles. Before finishing this remark, we only want recall a last B\"ar relation \cite[Lemma 2.1]{Ba2}:
\begin{equation}\label{GaussD}
{\bf D}\psi= \sum_{j=1}^{n}\gamma^{{\bf S}}(u_j)\nabla^{{\bf S}}_{u_j}\psi+\gamma^{\bf S}(N)\nabla^{\bf S}_N\psi=
\frac{n}{2}\gamma(\overrightarrow{H})\psi-\sum_{j=1}^{n}\gamma(u_j)\nabla_{u_j}\psi,
\end{equation}
where $\overrightarrow{H}$ is the mean curvature  vector field defined on $\Sigma$  and $\psi\in\Gamma({\mathbb S}\Omega_{|\Sigma})$ and $\{u_1,\dots,u_n\}$ a positive local orthormal frame of $\Sigma$. This equality  is a vectorial version of \eqref{bD}.  
}
\end{remark}

\section{A spinorial Reilly inequality}
A basic tool to relate the Dirac operator with the geometries of a connected compact spinorial manifold $\Omega$ and  those of their boundaries $\Sigma$, not necessarily connected,  that will be use, as in the closed case (see \cite{Fr}), is the integral version of the aforementioned Schr{\"o}dinger-Lichnerowicz formula \eqref{Dq}.
 In fact, from definitions above, one obtains the integral formula
 \begin{equation}\label{IR}
\int_{\Sigma}\left( ({\bf D}\psi,\psi)-\frac{n}{2}H|\psi|^2\right)=\int_\Omega\left(|\nabla\psi|^2-|D\psi|^2+\frac{1}{4}R|\psi|^2\right).
\end{equation}
In this situation, we will consider the pointwise spinorial Schwarz inequality$$
|D\psi|^2\le (n+1)|\nabla\psi|^2,\qquad\forall \psi\in\Gamma(\S\Omega),$$
where the equality is achieved only by the so-called {\em twistor spinors} (see \cite{BHMM,BFGK,Fr,LM,Tr}). Then we get the following integral inequality, called {\em Reilly inequality} (see \cite{ HMZ1}, for example)  because of its similarity with the corresponding one  obtained in \cite{Re} for the Laplace operator. Namely
\begin{equation}\label{Reilly}
\int_{\Sigma}\left( \left<{\bf D}\psi,\psi\right>-\frac{n}{2}H|\psi|^2\right)\ge\int_\Omega\left(\frac{1}{4}R|\psi|^2-\frac{n}{n+1}|D\psi|^2\right),
\end{equation}
where the equality is reached only for twistor spinors $\psi$ on $\Omega$ . The structure of this inequality gives rise to the consideration of boundary problems for the first order elliptic Dirac operator $D$ of the bulk manifold $\Omega$. To do this, we need suitable (elliptic) boundary conditions which cannot be   of Dirichlet or Neumann type. Think, for example, in the problem of finding a holomorphic function on a unit planar disc by prescribing its restriction to the unit boundary circle. The series expansion of the holomorphic function allows to prescribe exclusively {\em half} of the Fourier coefficients of the restriction function to the boundary. In fact, the elliptic boundary conditions for $D$ should  to determine one suitable {\em half} of  the values of the spinors fields  restricted to the boundary.

\section{The APS boundary condition}
Let us say some words about the {\em well-posed} boundary problems associated to the Dirac operator $D$ on an $(n+1)$-dimensional   compact connected spinorial Riemannian $\Omega$ with non-empty (not necessarily connected) $n$-dimensional boundary $\Sigma$. In this case, the boundary $\Sigma$ is our (embedded) hypersurface. It is a  well-known fact that the Dirac operator $D:\Gamma(\S\Omega)\rightarrow\Gamma(\S\Omega)$ on a compact spinorial Riemannian  manifold $\Omega$ with non-empty boundary $\Sigma$  has an infinite dimensional kernel and a closed  image with finite codimension. It is fundamental lo look for conditions $B$ to be imposed on the restrictions to the boundary $\Sigma$ of the spinor fields on $\Omega$ so that this kernel becomes finite-dimensional and, then, the problem 
\begin{equation}\label{boun-pro}\tag{BP}
\left\{
\ba{lll}
D\psi&=\Phi&\qquad\hbox{on $\Omega$}\\
B^{APS}\psi_{|\Sigma}&=\chi&\quad\hbox{along $\Sigma$},
\ea \right. 
\end{equation}
 is of Fredholm type. In this case, we will have  smooth solutions for any data $\Phi$ and for $\chi$ belonging to a certain subspace with finite codimension. These solutions will be unique up to the finitude of the dimension of the kernel. 

To our knowledge, the study of elliptic boundary conditions suitable for an elliptic operator $D$ (of any order, although for simplicity, we only consider first order operators) acting on smooth sections of a Hermitian vectorial bundle $F\rightarrow \Omega$ began in the fifties of past century by Lopatinsky and Shapiro (\cite{Ho,Lo}), but the main tool was discovered by Calder{\'o}n in the sixties: the so-called {\em Calder{\'o}n projector}.
This is a pseudo-differential operator of order zero  with principal symbol $\calp_+(D):T\Sigma\rightarrow {\rm End}_\C(F)$ which depends only on the principal symbol $\sigma_D$ of the operator $D$ (see  \cite{BW,Se} for details). One of the more important features of the Calder{\'o}n projector is that its principal symbol detects the  {\em ellipticity of a boundary condition}, or in other words, if the corresponding boundary problem \eqref{boun-pro} is a {\em well-posed problem} (according to Seeley in \cite{Se}). The reader interested in more details about the of boundary criteria for the ellipticity could see \cite{Se} or \cite[Chap. 18]{BW}. 
Most of the proofs of results in these two works implicitly suppose that the metric of $\Omega$ is a product in a collar neighbourhood of the boundary $\Sigma$ (it is clear in \cite{BW}), although the results remain to be true without this assumption.  A more modern and clearer study about the ellipticity of boundary problems for the Dirac operator, where we easily see that the metric of $\Omega$ close to $\Sigma$ does not have to be a product, can be found in \cite{BaBa1,BaBa2}.

Among the different global and elliptic local boundary conditions for the Dirac operator, the first  to be discovered and used, and the best known is the so-called Atiyah, Patodi and Singer condition. It was introduced in \cite{APS}  in order to establish the Index Theorem for compact manifolds with non-empty
boundary. Later, this condition has been revealed useful to study the positive mass
and the Penrose inequalities (see \cite{He2,Wi}). 

Precisely, this condition can be described as follows. We define the linear operator $B^{APS}:L^2({\bf S})\rightarrow
L^2({\bf S})$ as the orthogonal projection onto the subspace spanned by the non-negative eigenvalues of the self-adjoint intrinsic operator ${\bf D}$. Atiyah, Patodi and Singer showed in \cite{APS} (see also \cite[Prop.\!\! 14.2]{BW}) that $B^{APS}$ is a zero order pseudo-differential operator satisfying  the Lopatinsky-Shapiro  conditions required for the ellipticity to be satisfied.

Once the ellipticity of the APS boundary condition was established, it was proved the  fact (see \cite{APS,Se,BW,HMRo,BaBa1,BaBa2}) that  this spectral projection $B^{APS}$ from ${\bf D}$ provided an elliptic  boundary condition on $\Omega$ for $D$.
Then, we can prove that the homogeneous problem associated to \eqref{boun-pro}  admits only the trivial solution if the scalar curvature of $\Omega$ satisfies $R\ge 0$ and the mean curvature $H$ of $\Sigma$ is also non-negative. So \eqref{boun-pro} is Fredholm with unique solutions for suitable data. Moreover, one can check that, if we add the assumption that $\Sigma$ either supports no harmonic spinors or $H>0$, any set of data $\Phi$ and $\chi$ in $\Gamma({\bf S})$ are admissible. The last step is to including all this information in the Reilly spinorial inequality \eqref{Reilly}. One can see the details in the correct part of the failed \cite{Mon2}. 
 
\section{Two spinor integral inequalities obtained by using the APS boundary condition} 
Now, we will use the above spinorial machinery, and precisely some solutions of \eqref{boun-pro} subjected to APS boundary condition, to prove two results ( \cite[Theorem 2, Theorem 7]{Mon2}),  focus on elaborating a proof of the positivity theorem by Shi and Tam \cite{ST1,ST2} not  having recourse to glue infinite pieces to $\Omega$ or to solve difficult non-linear equations or taking for granted the positivity of the ADM mass demonstrated  by Schoen and Yau \cite{SY}  and Witten \cite{Wi}. We will get a spinorial proof of the positivity for the Brown-York mass that it is possible to generalise to domains enclosed by {\em mean-convex} boundaries and, as a consequence, a rigidity theorem for this type of Euclidean bodies.

 We will apply {\em \`a la Reilly} these results of  spinorial nature, thought our interest, more than the resolution of \eqref{boun-pro}, is  solving these equations on the bulk manifold for  using them to obtain results on its boundary.

\begin{th1}
Let $\Omega$ be a compact connected spinorial Riemannian manifold of dimension $n+1$ with non-negative scalar curvature $R\ge 0$ and having a non-empty boundary $\Sigma$ whose inner mean curvature $H\ge 0$ is also non-negative (mean-convex). Suppose also that $\Sigma$ either does not support harmonic spinors or $H>0$. Then, for every spinor $\phi\in H^1({\mathbb S}\Sigma)$, we have$$
\int_\Sigma|{\bf D}\phi|\,|\phi| \ge \frac{n}{2}\int_\Sigma H|\,\phi|^2.$$
The equality occurs if and only if    $\phi$ is the restriction to $\Sigma $ of a parallel spinor on $\Omega$ and $H$ is constant. Hence, using \eqref{nablaS}, 
we have
$$\nabla^{\bf S}_X\phi=-\frac1{2}\gamma^{\bf S}(AX)\phi\quad\hbox{\rm and}\quad {\bf D}\phi=\frac{n}{2}H\phi,$$
for every vector field $X$ tangent to $\Sigma$. So, the eigenspace associated to $\frac{n}{2}H$ consists in the restriction to $\Sigma$ of the space $E(\frac{n}{2})$ of parallel spinors defined on $\Omega$ (and the one associated to $-\frac{n}{2}$ will be $\gamma(N)\big(E(\frac{n}{2})\big)$.
\end{th1}

\pf

Let take us  a spinor $\phi$ defined on $\Sigma$ and $\xi$ the solution to \eqref{boun-pro} with data $\Phi=0$ and $B^{APS}\xi=B^{APS}\phi$. The hypotheses $R\ge 0$ and $\Phi=0$ imply that the right hand side of the Reilly inequality \eqref{Reilly} is non-negative. Then, using the pointwise Schwart inequality, we have 
$$
\int_\Sigma|{\bf D}\phi|\,|\phi| \ge \frac{n}{2}\int_\Sigma H|\,\phi|^2,
$$
as we were trying to prove.
In the case of equality, both sides of the integral equality \eqref{IR} vanish. Then $\Phi$ is harmonic and twistor-spinor on $\Omega$. Hence, it is a parallel spinor \cite{W,Ba1,Gi} and it is an open piece on the some of complete Ricci-flat manifolds listed   by Wang in \cite{W}. Moreover, from \eqref{bD}, we have the two other assertions in the statement on the theorem. Lasty, the vanishing of the left hand side implies that $H$ is constant because the equality is attained in the Schwarz inequality. Then $\phi$ is an eigenspinor of ${\bf D}$ for the eigenvalue $\frac{n}{2}H$ and the corresponding eigenspace consists in the restrictions to $\Sigma$ of the space of parallel spinors defined on $\Omega$. 
Then, if the boundary $\Sigma$ was not connected we would choose $\xi$ such that $\xi=0$ along of one of the connected components. Then, as $\xi$ is parallel, we would $\xi\equiv 0$.
\hfill $\square$

\begin{th2}
Let $\Omega $ be an $(n+1)$-dimensional compact connected spinorial Riemannian manifold whose scalar curvature satisfies $R \ge n(n+1)$  
and having a non-empty boundary $\Sigma  $ whose inner mean curvature $H \ge 0$ is  non-negative (mean-convex).  Then, for every spinor $\phi\in H^1({\mathbb S}\Sigma)$, we have$$
\int_\Sigma\sqrt{|{\bf D}\phi|^2\,|\phi|^2-\left<{\bf D}\phi,\gamma(N)\phi\right>^2} \ge\frac{n}{2}\int_\Sigma H\,|\phi|^2.$$
The equality holds if and only if the spinor field $\phi$
 is the restriction to $\Sigma$ of a $\pm \frac1{2}$-real Killing spinor $\psi$ defined on $\Omega$ (that is, $\nabla_X\psi=\pm\frac1{2}\gamma(X)\psi$ for all $X\in\Gamma(T\Omega))$. Thus, the scalar curvature $R$ of $\Omega$ must be identically $n(n+1)$, $H$ is constant  and
 $$
 {\bf D}\phi=\frac{n}{2}H\phi\pm\frac{n}{2}\gamma(N)\phi,
 $$
  In this case, 
 $$
|{\bf D}\phi|\ge|\left< {\bf D}\phi,\gamma(N)\phi\right>|=\frac{n}{2},
$$
with equality if and only if $H=0$.
Moreover,
$$
\xi^\pm=\left(H+\sqrt{1+H^2}\right)\phi \pm \gamma(N)\phi,
$$
 are eigenspinors of ${\bf D}$ associated to the eigenvalue $\frac{n}{2}\sqrt{1+H^2}$ (and, so, $\gamma(N)\xi^\pm$ are  eigenspinors associated to $-\frac{n}{2}\sqrt{1+H^2})$. That is, the eigenspaces associated to $\pm\frac{n}{2}\sqrt{1+H^2}$ are spanned by the real $\pm\frac1{2}$-Killing spinor on $\Omega$. 
\end{th2} 

\pf

In \cite{HMRo}, we, Hijazi and Rold\'an found a lower bound, in terms of $R$, for the largest non-positive and the smallest non-negative of the eigenvales for $D$ subjected to the boundary APS solution. Let us denote both by $\lambda_1(D)$. In fact, we saw that
$$
|\lambda_1(D)|^2>\frac{n+1}{4n}\min_\Omega R.
$$
Now, we make use of the Fredholm Alternative and the fact that $\Sigma$ has no harmonic spinors and, then, we deduce the existence of two spinors $\chi^\pm$ on $\Omega$ such that
 \begin{equation}\label{eigen-pro}\tag{EP}
\left\{
\ba{lll}
D\psi^\pm&=\mp\frac{n+1}{2}\psi^\pm&\qquad\hbox{on $\Omega$}\\
B^{APS}\psi_{|\Sigma}&=0&\qquad\hbox{along $\Sigma$}.
\ea
\right.
\end{equation} 
Putting this in the left hand side of the Reilly inequality \eqref{Reilly}, we deduce that the left hand side satisfies 
$$
0\le \int_{\Sigma}\left(| \left<{\bf D}\psi^\pm,\psi^\pm\right>|-\frac{n}{2}H|\psi^\pm|^2\right),$$
and the equality holds when $R=n(n+1)$ is constant and $\psi^\pm$ is a $|\frac1{2}|$-real Killing spinor.
Since the spinors $\psi^\pm$ and $\gamma(N)\psi^\pm$ form an orthonormal pair, we have
$$
| {\bf D}\psi^\pm|^2|\psi^\pm|^2\ge  \left<{\bf D}\psi^\pm,\psi^\pm\right>^2+\left<{\bf D}\psi^\pm,\gamma(N)\psi^\pm\right>^2.$$
From the two last inequalities, we obtain the integral inequality of the statement of the theorem.

In the case of the equality, $\psi^\pm$ must be an eigenspinor for $D$ associated to $\mp\frac{n+1}{2}$ and, furthermore, a twistor-spinor. Then, it is a $\pm\frac1{2}$-Killing spinor and, from \eqref{bD},
$$
 {\bf D}\phi=\frac{n}{2}H\phi\pm\frac{n}{2}\gamma(N)\phi.
 $$
Since the equality is also achieved in the pointwise Schwarz inequality, $H$ has to be constant and $|{\bf D}\xi|^2=
(1+H^2)|\xi|^2$ and, so, there are no harmonic spinors on $\Sigma$. The remaining assertions in the statement in the theorem are not too difficult to prove.\hfill $\square$

\begin{remark}{\em 
 It is worthy to note that both Theorems {\bf A} and {\bf B} above exhibe an almost total parallelism between them, except for a very precise point.  It is credible that first one should work 
for spinorial  Riemannian manifolds 
$\Omega$ whose tendency is to be Ricci-flat (to be endowed with a non-trivial parallel spinor field) and the second one for
 manifolds which tend to be Einstein manifolds with positive scalar curvature  $1$ (endowed with a real Killing spinor field). Recall that $n(n+1)$ is just the value 
of the scalar curvature of a unit sphere ${\mathbb S}^n$. This is why  the reader will have already noticed that, in  Theorem {\bf A}, the hypothesis about the 
scalar curvature on $\Omega$ is an inequality, namely $R\ge 0$ and, however, in Theorem {\bf B}, it is  another one $R\ge n(n+1)$.  

The relevant difference between the two theorems we talked about below refers to the hypothesis $H>0$ which. If we carefully examine the proofs of the two theorems, we will conclude that this assumption are not necessary in Theorem {\bf B}.

Of course, these inequalities are  the most popular  hypotheses in solved comparison theorems. We we will see that also will be able to improve proofs of already solved problems and to prove some unsolved conjectures. 

For example, among these latter 
and, on top of all, the inequality $R\ge n(n+1)$ is just the assumption imposed in the statement of well-known unsolved {\em Min-Oo's Conjecture}. This condition is necessary, because Huang and Wu proved in \cite{HW1} that, otherwise, one can 
perturb the hemisphere at an interior point so that $R\le n(n+1)-\varepsilon$, for some
small $\varepsilon > 0$, without changing the assumptions on the boundary \cite{HW1}. Moreover the original hypothesis  $R\ge n(n+1)$ has been invalidated as a sufficient hypothesis by the counterexamples built by Brendle, Marques and Neves \cite{BMN}. Furthermore, on their all counterexamples and the more recent built to \cite{Sw},   
it appears at least one point with this strict inequality in the bulk manifold. Hence,
we will bring to add some complementary hypothesis, namely, that the spin Riemannian structure on the boundary is the same that that one on equators of the ambient sphere.} It would interesting to check if the aforementioned counterexamples inherit from its seven-dimensional ambient Clifford multiplications different to the usual one. 
 \end{remark}
 
\section{The non-negative case}\label{NNC}
Now, we will take a  {\em more geometrical and less physical} look to the statements in Theorems {\bf A} and {\bf B}. Some of these results have been already published by ourselves and other authors 
and other ones are more or less slight improvements of known statements.   Also, we will obtain new geometrical results for surfaces immersed into spheres.
The first consequence of Theorem {\bf A} is almost immediate.

\begin{theorem}[\cite{HMZ1}, Theorem 6]\label{lambda1}
Suppose that $\Omega$ is a connected compact spinorial Riemannian manifold whose scalar curvature $R$ is non-negative with a non-empty (non-necessarily connected) boundary $\Sigma$ with mean curvature $H$ non-negative  (mean-convex). If either $\Sigma$ does not have harmonic spinors or $H>0$, then,$$
 |\lambda_1({\bf D})|\ge \frac{n}{2}\min_\Sigma H.$$
  If the equality holds, $\Omega$ is Ricci-flat and admits non-trivial parallel spinors. It is clear that $H=\frac{2}{n}\lambda_1({\bf D})$ is constant and the eigenspace corresponding to $\lambda_1({\bf D})$ consists of the restrictions to the boundary $\Sigma$ from non-trivial parallel spinor \cite{W,Ba1,Gi} defined on the bulk manifold. Moreover, $\Sigma$ has to be connected. Then, 
 (see \cite{W}), either $\Omega$ is flat and $N=\hbox{\rm mult\,}(\lambda_1({\bf D}))= 2^{[\frac{n+1}{2}]}$ for any $n\ge 1$, or it is a Calabi-Yau manifold 
 with $n\ge 4$ and $N=2$, or a hyperk\"ahler manifold with $n\ge 8$ and $N=\frac{n}{4}+1$, or a spin manifold such that $n=8$, $N=1$, supporting a 
  parallel three-form, or, finally,  a manifold with $n=7$, $N=1$,  endowed with a  parallel four-form.  \end{theorem}

\pf 

Let us pick an eigenspinor $\phi$ corresponding to $\lambda_1({\bf D})$. The hypotheses about $\Omega$ and $\Sigma$ permmet us to put it in the inequality of Theorem {\bf A}. Then
$$
\int_\Sigma| \lambda_1({\bf D})^2|\phi|^2 \ge \frac{n}{2}\int_\Sigma H|\,\phi|^2.
$$
From this, the lower estimate for $\lambda_1({\bf D})$ is a trivial consequence. In the equality is achieved, the same Theorem {\bf A} implies that $\Omega$ is a piece one of the simply connected Ricci-flat manifolds listed in \cite{W} and that eigenspace of  $\lambda_1({\bf D})$ consists of the restrictions to $\Sigma$ of the space of parallel spinors defined on $\Omega$.

 Finally, a few words about the connectedness of the boundary $\Sigma$. 
  Let $\phi$ be a non-trivial spinor field on $\Sigma$ with ${\bf D}\phi=
\frac{n}{2}H\phi$. Then, we may define another spinor field $\widetilde{\phi}$ by
$$
\widetilde{\phi}=
\left\{
\begin{array}{l}
\phi \quad\hbox{\rm on}\quad \Sigma_0\\
0\quad\hbox{\rm on}\quad \Sigma-\Sigma_0,
\end{array}
\right.
$$ 
where $\Sigma_0$ is a connected component of $\Sigma$. So, we get ${\bf D}\widetilde{\phi}=\frac{n}{2}H\widetilde{\phi}$, as well. Then $\widetilde{\phi}$ is a restriction of a parallel spinor field defined on $\Omega$. But, a non-trivial parallel spinor has non-null constant length. Hence, $\Sigma=\Sigma_0$.
\hfill $\square$ 

\begin{remark}\label{improve}{\rm
 If the bulk manifold satisfies $\ric= 0$ (for example, if 
 it supports a parallel spinor field), our lower estimate  improves, in the enclosing case, the well-known lower bound  by Fiedrich in \cite{Fr} 
 $$\lambda^2_1({\bf D})\ge \frac{n}{4(n-1)} \min_\Sigma R_\Sigma .$$
 The cause is the Gau{\ss} equation relating the scalar curvatures of $\Omega$ and its hypersurface $\Sigma$
 $$
0=   R-2\,\ric(N,N)=R_\Sigma-n^2H^2+|A|^2\ge R_\Sigma-n(n-1)H^2 ,
 $$
 where $R_{\Sigma}$ and $A$ are the scalar curvature function of $\Sigma$ and the second fundamental form with respect to $N$, respectively. Then
 $$
 H^2\ge \frac1{n(n-1)}R_\Sigma 
 $$
 and, so,
 $$
 \frac {n^2}{4}H^2\ge \frac{n}{4(n-1)}R_\Sigma.
 $$}
\end{remark} 

\begin{proposition}\cite[Abstract and Theorem 3]{HM}\label{Habajo}
Suppose that $\Omega$ is a compact connected spinorial Riemannian  $(n+1)$-manifold with $R\ge 0$ and a non-empty (not necessarily connected) boundary. Suppose that the boundary 
$\Sigma$ is strictly (inner) mean-convex, that is, $H>0$. Denote by  ${\bf D}$   the intrinsic Dirac operator of the boundary. Then 
$$
\int_\Sigma\left(\frac{|{\bf D}\phi|^2}{H}-\frac{n^2}{4}H|\phi|^2\right) d\Sigma\ge 0,$$
for all spinor field $\phi$ on $\Sigma$.   Equality occurs in exactly the same cases as in Theorem \ref{lambda1} above (so, the boundary $\Sigma$ is connected). 
 \end{proposition}
\pf \;  Let us put 
$$
|{\bf D}\phi||\phi|=\frac{|{\bf D}\phi|}{\sqrt{H}}\sqrt{H}|\phi|,
$$
for any spinor field on $\Sigma$. Using the $L^2$-Schwartz inequality in  $\Sigma$, we have
\begin{equation}\label{Sch}
\int_\Sigma |{\bf D}\phi||\phi|\le \left(\int_\Sigma \frac{|{\bf D}\phi|^2}{H}\right)^{\frac1{2}}
\left(\int_\Sigma H|\phi|^2\right)^{\frac1{2}}
\end{equation}
and, of course, the equality holds only when the functions $|{\bf D}\phi|$ and $H|\phi|$ are colinear. Using \eqref{Sch} in Theorem {\bf A}, we have
$$
\left(\int_\Sigma \frac{|{\bf D}\phi|^2}{H}\right)^{\frac1{2}}
\left(\int_\Sigma H|\phi|^2\right)^{\frac1{2}}\ge\int_\Sigma|{\bf D}\phi|\,|\phi| \ge \frac{n}{2}\int_\Sigma H|\,\phi|^2.
$$
Simplifying, we get the  integral inequality searched for:
$$
\int_\Sigma \frac{|{\bf D}\phi|^2}{H}
\ge \frac{n^2}{4}\int_\Sigma H|\,\phi|^2.
$$
The equality is attained like in Theorem {\bf A} and, in this case, $\Sigma$ must be connected and have 
constant mean curvature $H=\frac{2}{n}\lambda_1({\bf D})$. \hfill $\square$

\begin{remark}\label{confor}{\rm We know that, on a spinorial Riemannian manifold $\Omega$, given two spin structures corresponding to two conformal  
Riemannian metrics $g$ and ${\bar g}$ on it, namely ${\bar g}=e^{2u}g$, the associated spinor bundles are identifiables. Then, if ${\bar {\bf D}}$ and 
${\bf D}$ are  Dirac operators on each one of them, we have the following conformal pseudo-invariance:
\begin{equation}\label{confD}
{\bar{\bf D}}(e^{-\frac{n}{2}u}\phi)=e^{-\frac{n+2}{2}u}{\bf D}(\phi),
\end{equation}
for all spinor field $\phi$ on the $(n+1)$-dimensional manifold $\Omega$. It was discovered by Hitchin in \cite[Section I.6]{Ht}. Then, using \eqref{confD}, one can convince 
him/her-self that the statement in Proposition \ref{Habajo} is equivalent to affirm that
$$
\lambda_1({\bf D}_{H^{-2}g})\ge\frac{n}{2}
$$
(remind that $H>0$) and that the equality occurs in the very same cases.
}
\end{remark}

 Now, suppose that the $(n+1)$-dimensional manifold $\Omega$ is endowed with a parallel (unit) spinor field $\psi_0$ (so, it is Ricci-flat and makes part of the list of manifolds in \cite{W} or in Theorem {\bf A}). Let us pick a vector field $X$ tangent to the manifold. With these two ingredients, we may define a new spinor field by means of the relation 
$$
\psi=\gamma(X)\psi_0.
$$ 
Then, it is immediate that
$$
\nabla_u\psi=\gamma(\nabla_uX)\psi_0,
$$
 Applying \eqref{bD}, we obtain
 $$
 {\bf D}\psi=\frac{n}{2}H\psi+({\div_\Sigma X})\gamma(N)\psi_0.
 $$
 Let us compute the norms of the right hand side on equation above. Thus
 \begin{eqnarray*}
|{\bf D}\psi|^2&=&\frac{n^2}{4}H^2|\psi|^2+({\div_\Sigma X})^2+nH({\div_\Sigma X})\left<\psi,\gamma(N)\psi_0\right>\\
&=& \frac{n^2}{4}H^2|\psi|^2+({\div_\Sigma X})^2+nH({\div_\Sigma X})\left<\gamma(N)\psi_0,  \gamma(X)\psi_0\right>\\
&=& \frac{n^2}{4}H^2|\psi|^2+({\div_\Sigma X})^2+nH({\div_\Sigma X})\left<N, X\right>.
\end{eqnarray*}
  From this, reminding that we are assuming that $H>0$, we have
  $$
\frac{ |{\bf D}\psi|^2}{H}=\frac{n^2}{4}H|\psi|^2+\frac{({\div_\Sigma X})^2}{H}+n({\div_\Sigma X})\left<N, X\right>.  
$$
By putting this in Proposition \ref{Habajo}, we obtain the following inequality
 
 \begin{equation}\label{ineq}
 0\le  \int_\Sigma\left(\frac{(\hbox{\rm  div$_\Sigma$}X)^2}{H}+n(\hbox{\rm  div$_\Sigma$}X) \left<N,X\right>\right)  
\end{equation}
for  all vector field $X$ tangent to $\Omega$. 

We will pick an arbitrary smooth function $f\in C^\infty (\Omega)$ defined on the bulk. Now, let us take $X$ as the gradient $X=\nabla g$ of the solution to the 
solution to the Neumann boundary equation
$$
\left\{
\begin{array}{l}
\Delta g={f},\quad\hbox{\rm on $\Omega$}\\ \\{\displaystyle
\frac{\partial g}{\partial N}=-\frac1{A(\Sigma)}\int_\Omega f, \quad\hbox{\rm along $\Sigma.$}}
\end{array}
\right.
$$
Then, 
\begin{equation}\label{CF}
\div X=\Delta g={f}=\div_\Sigma X+\left< \nabla_NX,N\right>.
\end{equation}
Let us add a new assumption on the arbitrary vector field $X$. {\em From now on, we will suppose that $X$ is a conformal field, that is, a flow of infinitesimal conformal
transformations}. This means that $X$ behaves in this way:
$$
L_X\left<\,,\right>=\frac1{n+1}({\div X})\left<\,,\right>=\frac1{n+1}f\left<\,,\right>,
$$ 
where $\left<,\,\right>$ is the Riemannian metric of $\Omega$. It is worthy to remark that two-dimensional manifolds are the only ones with an infinite-dimensional conformal group. In general, conformal fields have divergence with constant Hessian. But, this is a more or less known fact.  Moreover, it can be seen there that this constant is zero only if the field is parallel. Hence, 
$$
\left<\nabla_NX,N\right>=(L_X\left<\,,\right>)(N,N)=\frac1{n+1}({\div X})=\frac1{n+1}f.
$$
Then, when $X$ is conformal, the equation \eqref{CF} can be rewritten as follows
$$
\div X=\Delta g={f}=\div_\Sigma X+\frac1{n+1}f.
$$
Hence
\begin{equation}\label{divS}
\div_\Sigma X=\frac{n}{n+1}\div X=\frac{n}{n+1} f.
\end{equation}
the inequality \eqref{ineq} becomes 
\begin{equation}\label{ineq2}
0\le  \int_\Sigma\left(\frac{f^2}{(n+1)H}+ f \left<N,X\right>\right). 
\end{equation}
But, by the Divergence Theorem again,
\begin{eqnarray*}
\int_\Sigma f \left<N,X\right>&=&\int_\Sigma \left<N,fX\right>=-\int_\Omega \div (fX)\\
&=&-\int_\Omega f^2-\int_\Omega X\cdot\div X\\
&=&-\int_\Omega f^2-\int_\Omega \left< \nabla f,\nabla g\right>\\
&=&-\int_\Omega f^2+\int_\Omega  f \Delta g+\int_\Sigma f \frac{\partial g}{\partial N}.
\end{eqnarray*}
As $\Delta g=f$, we have finally
$$
\int_\Sigma f \left<N,X\right>=-\frac1{A(\Sigma)}\left(\int_\Sigma f\right) \left(\int_\Omega f\right).$$
We put this information in the inequality \eqref{ineq2} and have
$$
0\le  \int_\Sigma\frac{f^2}{H}-\frac{n+1}{A(\Sigma)}\left(\int_\Sigma f \right)\left(\int_\Omega f\right).  
$$

All these calculations and comments above constitute  the evidence for a type of Heintze-Karcher Inequality, so called because of 
 its similitude with the inequalities gotten in \cite{HK} (cf.\! \cite[Theorem 6.16]{MR3} and \cite[(4), p.\;251]{Br1}). It was Desmonts in \cite{De} the first author 
 which introduced the spinorial tools to study this type on inequalities (see also \cite{Re} and \cite[Chapter 6]{MR3}.
 
\begin{theorem}[Heintze-Karcher-Ros Inequality, \cite{HK,Ro}]\label{HK}
Suppose that $\Omega$ is a compact connected spinorial Riemannian  $(n+1)$-manifold with $R\ge 0$ and a non-empty (not necessarily connected) strictly mean-convex boundary. Suppose also that $\Omega$ is endowed wit a parallel spinor and a  conformal vector field $X$. If the boundary 
$\Sigma$ is, in strictly (inner) mean-convex, that is, $H>0$, then
$$
0\le  \int_\Sigma\frac{(\div X )^2}{H}-\frac{n+1}{A(\Sigma)}\left(\int_\Sigma \div X \right)\left(\int_\Omega \div X\right).
$$
The equality is attained like in Proposition \ref{Habajo} or Theorem {\bf A}. 
\end{theorem} \vspace{-.7cm}\hfill$\square$

\begin{remark}{\rm 
In general, the divergence of a conformal field has vanishing Hessian (it is an affine function), but if the manifold is Ricci-flat,  we have that $\div X$ is a non-null constant. In this case, the inequality above takes the form
$$
\hbox{\rm vol}(\Omega)\le \frac1{n+1}\int_\Sigma \frac 1{H}.
$$   
The equality is reached if and only if $\Sigma$ is an embedded connected and totally umbilical  into the manifold $\Omega$.
}
\end{remark}

As an immediate consequence, we obtain a celebrated result in Theory of Submanifolds. One can read the original proof, based on a maximum principle 
for non-linear PDE's, written by the author in \cite{Al}, we advise the lector to see it in \cite{Sp}. Also, it could be useful \cite{MR2} for a proof without spinors and generalisations to other ambient spaces.

\begin{theorem}[Alexandrov Theorem, \cite{Al,Ro,MR2,HMZ1}]\label{AT}
Suppose that $\Sigma$ is an enclosing hypersurface in a compact connected spinorial Riemannian  $(n+1)$-manifold endowed with a non-trivial parallel spinor (for example, flat spaces and cones on Sasaki, Einstein-Sasaki and 3-Sasaki manifolds, nearly-K\"ahler or seven-dimensional manifolds with a three-form $\omega$ with $\nabla\omega=\ast\omega$, (see \cite{Ba1}).  Then, if the mean curvature $H$ of $\Sigma$ is constant, then $\Sigma$ is  a sphere embedded in a totally umbilical way and lying in a fundamental domain of the universal cover of the manifold $\Omega$. \end{theorem}
\pf \!
It suffices to accept that, if $H$ is constant, then $H>0$ for the inner orientation on $\Sigma$, because $\Sigma$ lies we in a Ricci-flat manifold. Then, we can write Theorem \ref{HK}  adding this information. In fact, if $H>0$ is constant, Theorem \ref{HK} says
$$
(n+1)H\hbox{\rm vol\,}(\Omega)\le A(\Sigma).
$$
and the equality occurs like in Theorem \ref{HK}. On the other hand, \eqref{divS} implies
$$
\div_\Sigma X^\top=\div_\Sigma\left(X-\left<X,N\right>N\right)=\frac{n}{n+1}f+nH\left<X,N\right>.
$$
Integrating this equality and using the Divergence Theorem, we obtain
$$
0=\frac1{n+1}\int_\Sigma f+H\int_\Sigma\left<X,N\right>=
\frac1{n+1}\int_\Sigma f-H\int_\Omega f.
$$
Since, as in Theorem \ref{HK} above, we know that $f=\div X$ is a non-null constant,  we have really
$$
0=\frac {A(\Sigma)}{n+1}-H\hbox{\rm vol\,}(\Omega).
$$
Thus, we always have the equality in the inequality quoted above. All the remaining affirmations relative to the equality have already analised or one can find them in \cite{Mon1}. \hfill $\square$   

\begin{remark}{\rm 
It is also worthy to note that, as a consequence of the estimate in Theorem \ref{lambda1}, if $\Sigma$ admits an isometric and isospin immersion in an Euclidean space ${\mathbb R}^{n+k}$, $k\ge 0$, we can enhance our result. In fact, in this case, there exists an almost evident {\em upper} bound for $\lambda_1({\bf D})$. The proof simply consists of computing the Rayleigh quotient for the spinor field $\gamma(p)\phi_0$, where $p$ is the position vector and $\phi_0$ is a non-trivial constant vectorial function taking values in the ambient Euclidean space. In fact, from Remark \eqref{highD}, we know (see \cite[Section 5.2]{Gi}) that
\begin{equation}\label{Hp}
|\lambda_1({\bf D})|^2\le \frac{n^2}{4}\frac{{\displaystyle\int_\Sigma |{\overrightarrow H}|^2}}{\hbox{\rm vol}(\Sigma)}
\end{equation}
and the equality occurs if and only if $\phi$ is an eigenfield for $D$ associated to the eigenvalue $\frac{n}{2}$ and $|{\overrightarrow H}|=1$. This is equivalent to
$$
D\phi=\frac{n}{2}\gamma({\overrightarrow H})\,\phi+n\,\gamma(p)\phi_0=\frac{n}{2}\,\phi.
$$ 
Taking norms in this equality, we see that $\big<{\overrightarrow H},p\big>=-1$. Since ${\overrightarrow H}$ and $p$ are unit vectors, we have ${\overrightarrow H}=-p$ on $\Sigma$.  By taking derivatives in direction $u\in T\Sigma$, we have
$$
\nabla^\perp {\overrightarrow H}=0\quad\hbox{\rm and}\quad A_{\overrightarrow H}=I.
$$
So, one comes to the conclusion that the equality \eqref{Hp} above attains the equality if and only if $\Sigma $ is immersed  into a Euclidean space ${\mathbb R}^{n+\ell}$,
$1\le \ell\le k-1$, with non-trivial {\em parallel mean curvature} (cf.\! \cite[pp.\! 90-91]{Gi}). There are a plenty of papers by many authors searching about this kind on submanifolds. Among them, the works by Chen \cite{Ch} and, independently, by Yau \cite{Y1} allow us to assert that either $\Sigma$ in a minimal hypersurface in ${\mathbb S}^{n+k-1}$ or $n=2$ and it is either a unit two sphere in ${\mathbb R}^3\subset {\mathbb R}^m$ or, finally, a surface with constant mean curvature in  ${\mathbb S}^3\subset {\mathbb R}^4\subset {\mathbb R}^m$. When $\ell=1$, that is, $\Sigma$ is a hypersurface of ${\mathbb R}^{n+1}$, the two equations above can be easily understood. Indeed, they are equivalent to $H=1$, where now $H$ is the mean curvature {\em function} on $\Sigma$. Then, the Alexandrov Theorem \ref{AT} implies that $\Sigma$ is a unit $n$-sphere.}
\end{remark}

The following result is another consequence of our lower estimate for the first eigenvalue of the Dirac operator stated in Theorem \ref{lambda1}. In fact, in its original form, our assertion will be an adaption to the flat realm to an old conjecture proposed by Min-Oo in the spherical ambient (see \cite{Br3} and Theorem \ref{FMOo} below).  

\begin{theorem}[Flat Min-Oo's Conjecture]\label{FMOo}
Suppose that $\Sigma$ is an enclosing boundary hypersurface in an  $(n+1)$-dimensional connected compact spinorial manifold $\Omega$ with non-negative scalar curvature $R\ge 0$. If $\Sigma$ admits an isometric and isospin immersion into a Euclidean space ${\mathbb R}^{n+k}$, with $k\ge 1$   such that its mean curvature  accomplishes $|{\overrightarrow {H_0}}|\ge 1$, then $\Omega$ 
must be  a disc ${\mathbb D}^{n+1}\subset {\mathbb R}^{n+1}\subset {\mathbb R}^{n+k}$  and $\Sigma$ is the round sphere ${\mathbb S}^n$ enclosing it. 
\end{theorem}
\pf\! 
Combining Theorem \ref{lambda1} and Remark \ref{highD} with \eqref{Hp} and using that $|{\overrightarrow {H_0}}|\le 1$, we have 
$$
\frac{n^2}{4}\le \frac{n^2}{4}\min_\Sigma |{\overrightarrow H}|^2\le |\lambda_1(D)|^2\le \frac{n^2}{4}\frac{{\displaystyle\int_\Sigma |{\overrightarrow {H_0}}|^2}}{\hbox{\rm vol}(\Sigma)}\le \frac{n^2}{4}.
$$
Once attained the equality, on the one hand, there a parallel spinor  defined on $\Omega$ and this compact domain lies in one the five types of spin manifolds listed in Theorem \ref{lambda1}. Moreover, according to the digression before the statement of the theorem, we have
$$
D\phi=\frac{n}{2}H\gamma(N)\phi=\frac{n+1}{2}\gamma({\overrightarrow H_0})\phi.
$$
A first consequence is the relation $|{\overrightarrow H_0}|=H=1$. So, $H>0$ and we can apply Theorem {\bf A}. Furthermore, on one hand, as the  eigenspace $\lambda_1(D)$ consists of the restrictions to $\Sigma$ of the space of parallel spinors on $\Omega$, then its dimension is $N=2^{\left[\frac{n+k}{2}\right]}, 2, \frac{n}{4}+1, 1$, according  to $\Omega$ is flat, Calabi-Yau, hyperk\"aler or, finally,  one of the $7$ and $8$-manifolds with special holonomy group.  On the other hand, this same eigenspace is $N=2^{\left[\frac{n+k}{2}\right]}$-dimensional, because it comes from the space of parallel spinor of the ambient ${\mathbb R}^{n+k}$. After an easy checking, all the remaining cases can be reduced to three alternatives: a) 
$(n+1,k)=(\hbox{\rm arbitrary},1)$, $(n+1,k)=(\hbox{\rm arbitrary},0)$; b) $(n+1,k)=(2,0)$ or $(n+1,k)=(1,1)$   ; c) $(n+1,k)=(4,1)$. But all of them are impossible because of dimensional reasons (revisit the list in Theorem \
\ref{lambda1}). Thus, lastly, the unique remaining possibilities are that either $\Omega$ be a unit disc in ${\mathbb R}^{n+1}$ enclosed by $\Sigma={\mathbb S}^n$ or $\Omega$ is the domain enclosed in ${\mathbb S}^{n+1}\subset {\mathbb R}^{n+2}$ by any embedded minimal hypersurface $\Sigma$.  But, in this later case, $\Omega$ is a spherical domain and it cannot be endowed with on parallel spinors.\hfill $\square$

\begin{corollary}[``Original'' Flat Min-Oo's Conjecture]
Suppose that $\Sigma$ is an enclosing  hypersurface of an  $(n+1)$-dimensional Euclidean space ${\mathbb R}^{n+1}$. Moreover, it is  isometric and isospin to a unit sphere ${\mathbb S}^n$ and its  inner mean curvature $H$ satisfies $H\ge 1$. Then $\Omega$ must a round disc ${\mathbb D}^{n+1}$ of radius $1$ and $\Sigma$ is the round sphere ${\mathbb S}^n$ enclosing it. 
\end{corollary}
\pf\;
 The  point here is to realise that the sphere ${\mathbb S}^n$ has an only {\em topological spin structure} which does not admit harmonic spinors and that $|\lambda_1^\pm|=\frac{n}{2}$ for its usual {\em Riemannian spin structure}. Finally, it is  evident that it admits an immersion in ${\mathbb R}^{n+1}$ with ${\overrightarrow H}=N$ (and, then, $H=1$). 
 \hfill $\square$

\begin{remark}{\rm
Our  answer to corresponding {\em flat} Min-Oo's conjecture was also independently obtained by Miao \cite{Mi1} (see Remark 1), although he makes no mention of spin structures, at all. Really, our  statement is also not the {\em original Min-Oo conjecture}. Indeed, it  is one of the possible adaptions of the well-posed conjecture to an Euclidean context. The original conjecture was posed for totally geodesic hypersurfaces of an  $(n+1)$-dimensional sphere. }
\end{remark}  

 \begin{remark}\label{metrispin}{\rm
 To understand that a spin Riemannian manifold can admit an only spin  {\em topological} structure and different {\em spin metric structures}, it suffices to think in the Euclidean space ${\mathbb R}^7$. In fact, since it is a
 Euclidean space, it has vanishing first and second Stiefel-Whitney, $w_1=0$ and $w_2=0$. Then, this seven-dimensional space, endowed as the usual Riemannian metric, is a spin manifold. Since $H_1({\mathbb R}^7,{\mathbb Z}_2)=0$, it has an only {\em topological} spin structure. However, it supports at least two {\em metrically} different ones, namely, the standard one, coming from the usual Clifford product and that induced when ${\mathbb R}^7$ is thougth of as the set of purely imaginary  numbers in the octonions space (see \cite{Ba1,ABF,W}). The first one has a maximal number $2^3=8$ of independent parallel spinors and, instead, the second one has only $1$ of that type of spinor fields \cite{Ba1}. Then, both spinorial structures are topologically equivalent, but not metrically.   
 }
 \end{remark}

In what follows, we will show the principal aim which led us to look for a  proof of positivity of the Brown-York, alternative to that of the pioneers Shi and Tam. We wanted a proof avoiding unbounded domains, complicated PDE's and the positivity of another mass: the so-called ADM-mass. This latter was thought for three-manifolds asymptotically flat (in a suitable and {\em ad hoc} sense) with non-negative scalar curvature. Also, we wondered if the convexity of the hypersurfaces was really necessary. It is from this quest that Theorems {\bf A} and {\bf B} grew out.  
It is obvious that the following result is a generalization of the positivity theorem for the Brown-York mass previously proved for strictly convex surfaces by Shi and Tam. In their proof, the solution of difficult boundary equations and the positivity of the ADM-mass obtained by Shoen-Yau and Witten \cite{SY,Wi}, in the context of asymptotically flat manifolds, are essential components. Here, these difficulties are avoided and, as we have already remarked somewhere above, this {\em compact version} of the theorem implies the asymptotically flat version for the ADM mass (see \cite{HMRa1}).

\begin{theorem}[Shi-Tam's Theorem: Brown-York mass for mean-convex surfaces]\label{ST} 
Let $\Omega$ be a compact spin Riemannian manifold of dimension $n+1$ with non-negative scalar curvature $R\ge 0$ and having a non-empty boundary $\Sigma$ whose inner mean curvature $H\ge 0 $ is  mean-convex and has no harmonic spinors. Suppose that there is an isometric and isospin immersion from $\Sigma$ into another spin manifold $\Omega_0$ endowed with a non-trivial parallel spinor field and let $H_0$ its mean curvature with respect to any of its orientations.  Then, we have$$
\int_\Sigma H\le \int_\Sigma |H_0|.$$
The equality implies that $H=|H_0|=H_0$. Then, if $n=2$, $\Omega_0$ is a domain in ${\mathbb R}^3$ and the two embeddings differ by a direct rigid motion
\end{theorem}

{\bf Proof.}  Denote by $\psi$ the parallel spinor on $\Omega_0$ and let $\phi=\psi_{|\Sigma}$ its restriction onto $\Sigma$ through the existent immersion. Let's recall that the parallelism of $\psi$ (see (22)) gives$$
{\bf D}\phi=\frac{n}{2} H_0\phi\qquad\hbox{and}\qquad |\phi|=1.$$  Suppose that $\Sigma$ has no harmonic spinors. Now, we apply Theorem 2 and have the desired inequality$$
\int_\Sigma H\le \int_\Sigma |H_0|.$$

If the equality is attained, then $$\frac{n}{2}H_0=D\phi=\frac{n}{2}H
$$ and so $H=H_0> 0$.
Then, the immersion of $\Sigma$ into the second ambient space $\Omega_0$ is strictly mean-convex as well (with respect to
the choice of inner normal to $\Omega$).  

When $n=2$, from this equality and the fact that $K=K_\phi$ (because the two embeddings are isometric and preserve the Gauss curvatures), we deduce that the two second fundamental forms coincide. The Fundamental Theorem of the Local Theory of Surfaces allows us to conclude that the two boundaries differ by a  rigid motion of the Euclidean space. \hfill$\square$
 
\begin{remark}{\em
Note that, in the original Shi-Tam original result, the authors assume that the boundary $\Sigma$ is strictly convex. Then, the well-known answers by Pogorelov \cite{Po} and Nirenberg \cite{Ni} to the Weyl problem guarantee the existence of a geometrically unique  isometric embedding into the Euclidean space as the boundary surface of a convex body. Instead, we need suppose the existence of this second isometric immersion with $H_0\ge 1$. Moreover, from the result {\em a fortiori} one deduces $H_0>0$} .
\end{remark}

\begin{corollary}[Alexandrov-Fenchel-Minkowski Theorem]\label{finta}
Let $\Omega$ be a compact Riemannian manifold of dimension three with non-negative scalar curvature $R\ge 0$ and having a mean-convex boundary $\Sigma$ isometric to a sphere of any radius.  Then, we have$$
\int_\Sigma H\le \sqrt{\pi A(\Sigma)}$$
where $A(\Sigma)$ is the area of $\Sigma$. It the equality holds, then the two boundaries are spheres of the same radius.
\end{corollary}
{\bf Proof.} First, observe that all Riemannian three-manifolds are spin. It is  clear that the boundary ${\mathbb S}^2$ of $\Omega$ admits and isometric and isospin (the sphere supports a unique spin structure, up to ${\mathbb S^6}$) embedding into the Euclidean space ${\mathbb R}^3$ with $|H_0|=1/r$ and area $A(\Sigma)=\pi r^2$, where $r>0$ is the radius of the sphere. The fact that the two embeddings of ${\mathbb S}^2$ are isometric allows us to finish. \hfill $\square$

\begin{remark}{\em 
The integral inequality in Corollary \ref{finta}, for strictly convex surfaces of ${\mathbb R}^3$,  is attributed to Minkowski (1901), although its very probable that it were previously known to Alexandrov and Fenchel. Recently, it has been proved  in \cite{DHMT}  that the Minkowski inequality is not valid for any compact surface, although they proved it is for the axisymmetric ones. It is also worthy to remark the following conjecture by Gromov: If $\Sigma$ is the boundary of a compact Riemannian manifold $\Omega$, then, if $R\ge \sigma$, for a certain constant $\sigma$, where $R$ is the scalar function of $\Omega$, then there exists a constant $\Lambda(\Sigma,\sigma)$ such that 
$$
\int_\Sigma H\le \Lambda(\Sigma,\sigma).$$}
\end{remark}

\begin{corollary}[Cohn-Vossen Rigidity Theorem for Mean-Convex Domains, \cite{Mon2}]\label{CV}
Two isometric and isospin strictly mean-convex compact surfaces or with $H=0$ in the Euclidean space ${\mathbb R}^3$ must be congruent.
\end{corollary}
{\bf Proof.} Let $\Omega$ and $\Omega_0$ be the two domains determined in ${\mathbb R}^3$ by two corresponding surfaces identified by means of an isometry. Then, we can apply Theorem \ref{ST} interchanging the roles of $\Omega$ and $\Omega_0$ and applying the case of the equality. \hfill $\square$

\begin{remark}{\em 
Much more recently, in the context \cite{SWWZ} of {\em fill in} problems posed  firstly by Bartnik, it has been proved that, if $\Omega$ is the hemisphere $B^{n+1}$ and $\gamma$ is a metric on the boundary ${\mathbb S}^n$ isotopic to the standard one with mean curvature $H>0$, then there is a constant $h_0=h_0(\gamma)$ such that$$
\int_\Sigma H\le h_0.$$ 
It is clear that this result and our Corollary \ref{finta}, together with the question by Gromov,  belong to a same family.} 
\end{remark}

\section{Ambients with positive scalar curvature}
Until now we have suppose that the scalar curvature of our compact spin Riemannian manifold $\Omega$ satisfied $R\ge0$ (Euclidean context). Let us enhance this positivity assumption to $R\ge n(n+1)$ (spherical context). This lower bound  is precisely the constant value of the scalar curvature of the $(n+1)$-dimensional unit sphere. 
Then, putting this assumption and the  Schwarz inequality
\begin{equation}
|D\psi|^2\le (n+1)|\nabla \psi|
\end{equation}
(already used in  Section 3)
into the right hand side of the 
Weitzenb\"ok-Lichnerowicz formula, we obtain Theorem {\bf B}.
\begin{equation}
\int_\Sigma\left(\left<{\bf D}\psi,\psi\right>
-\frac{n}{2}H|\psi|^2\right)\ge \int_\Omega\left(-\frac{1}{n+1}|{D}\phi|^2+\frac{n+1}{4}|\psi|^2\right),
\end{equation}
for all compact spin manifold $\Omega$, with equality  only for the twistor  spinor fields on $\Omega$. 
From this starting point,  by using this integral inequality, we will work in a similar, but a more geometrically and few more elaborated, way as in Theorem 2,  and will get an expanded version of  Theorem {\bf B}.

\begin{theorem}\label{lowerS}
 Let $\Omega $ be a $(n+1)$-dimensional compact spin Riemannian manifold whose scalar curvature satisfies $R \ge n(n+1)$  and having a non-empty boundary $\Sigma  $ without harmonic spinors. Then, for every spinor $\phi\in H^1(\Sigma)$, we have$$
\int_\Sigma\sqrt{|{\bf D}\phi|^2|\phi|^2-\left< {\bf D}\phi,\gamma(N)\phi\right>^2}\ge\frac{n}{2}\int_\Sigma H|\phi|^2 .$$
The equality holds if and only if $\phi$ is a spinor field coming from a $|\frac1{2}|$-real Killing  spinor $\psi$ defined on $\Omega$ (see \cite{Ba1} and \cite[Appendix A]{Gi}), and so
$$
\phi=\frac{n}{2}H\psi \pm\frac{n}{2}\gamma(N)\psi.
$$
Like in Theorem {\bf A}, the boundary $\Sigma$ must be connected and $H$ is constant in case of equality. The bulk $\Omega$ is contained either in an $(n+1)$-sphere with multiplicities $(2^{[\frac{n}{2}]},2^{[\frac{n}{2}]})$ for its $\pm$-Killing spinors, or in a $(4n+1)$-dimensional Einstein-Sasaki manifold, with $n\ge 1$ and multiplicities $(1,1)$, or in a $(4n+3)$-dimensional Einstein-Sasaki manifold but not 3-Sasaki, with $n\ge 2$ and multiplicities $(2,0)$,  or in a
 $(4n+3)$-dimensional 3-Sasaki manifold, with $n\ge 2$ and multiplicities $(n+2,0)$, or in a seven-dimensional manifold endowed with a 3-form $\omega$ such that $\nabla\omega=\ast\omega$ but not Sasaki with multiplicities $(1,0)$, or, finally, in a six-dimensional $(1,1)$-nearly-K\"ahler manifold with multiplicities $(1,1)$. \hfill$\square$ 
 \end{theorem}

As far as we know, from this Theorem \ref{lowerS}, we obtain a new accurate lower estimate for the first eigenvalue of the Dirac operator in a hypersurface lying in a context of scalar curvature positive (cf. \cite[Chapter 3]{Gi}). This lower bound was expected  once we  already found the corresponding ones for non-negative   (Theorem \ref{lambda1} above) and negative ambient spaces (see \cite{HMZ1} and \cite[Theorem 3.7.1]{Gi}).

\begin{theorem}\label{lowerSS}
 Consider a connected compact spinorial  Riemannian manifold $\Omega$ of dimension $n+1$ with scalar curvature $R\ge n(n+1)$ and mean-convex boundary. Suppose that $\Sigma$ either does not support harmonic spinors (for example, if $R_\Sigma>0$), and let $\phi$ be the eigenspinor corresponding to either the smallest positive eigenvalue or to the greatest negative one, namely  indistictly, $\lambda_1({\bf D})$,  of the intrinsic Dirac operator ${\bf D}$ of the boundary. A direct application of Theorem {\bf B} gives$$
\int_\Sigma\left(|\lambda_1({\bf D})|-\frac{n}{2}\sqrt{1+H^2}\right)|\phi|^2\ge 0.$$
As a consequence, we obtain the following lower bound:$$
 |\lambda_1({\bf D})|\ge \frac{n}{2}\min_\Sigma \sqrt{1+H^2}.$$
 This improves  the well-known intrinsic lower bound by Friedrich in \cite{Fr} for whatsoever spinorial Riemannian manifolds because the Gau{\ss} equation for the scalar curvatures in a unit $(n+1)$-dimensional manifold with $R=n(n+1)$  gives 
 $$
 \sqrt{\frac{nR_\Sigma}{4(n-1)}}\le \frac{n}{2}\sqrt{1+H^2}$$
 and the equality is attained only by the umbilical hypersurfaces. The equality holds if and only if $\phi$ is a spinor field coming from a  (positive or negative) real Killing  spinor $\psi$ defined on $\Omega$ (see \cite{Ba1} and \cite[Appendix A]{Gi}),
by
$$
\phi=\psi\mp \left(H\pm\sqrt{1+H^2}\right)\gamma(N)\psi,
$$
according to the sign of the Killing spinor on $\Omega$. \hfill $\square$
\end{theorem}

Surprisingly, from this theorem, we obtain a spinorial characterisation of the {\em Alexandrov embedded} minimal hypersurfaces in spin Riemannian manifolds with scalar curvature $R\ge n(n+1)$ (just the value which takes on the unit sphere ${\mathbb S}^{n+1}$). All of them have the same $\lambda_1({\bf D})$ and this value is peculiar uniquely to them. As a consequence, we can give a spinorial solution to the famous Yau's  \#100 Problem about minimal surfaces embedded in ${\mathbb S}^3$ (see \cite{Y2}) posed in 1982.

\begin{theorem}[Spinorial Version of Yau's \#100 Problem] \label{minimal}
 Let us consider a connected compact spinorial  Riemannian manifold $\Omega$ of dimension $n+1$ with scalar curvature $R\ge n(n+1)$ and mean-convex boundary. Suppose that $\Sigma$ does not support harmonic spinors and that $\lambda_1({\bf D})$  is either the smallest positive eigenvalue or to the greatest negative one  of the intrinsic Dirac operator ${\bf D}$ of the boundary,  both denoted by a same symbol. Then, we have the following    lower bound
 $$
 |\lambda_1({\bf D})|\ge \frac{n}{2},
 $$
and the equality is attained if and only if $\Sigma$ is minimal in $\Omega$. Moreover, in such a case, the corresponding eigenspace consists of the restrictions to $\Sigma$ of all the (positive or negative) real Killing spinor fields on $\Omega$. As for the multiplicity of $\lambda_1({\bf D})=\frac{n}{2}$, it suffices to invoke  Theorem {\bf B} above. 
\end{theorem} 
 \pf \! 
 It is clear that, under our hypotheses and Theorem \ref{lowerSS} above, we have
 $$
|\lambda_1({\bf D})|\ge \frac{n}{2}\max_\Sigma \sqrt{1+H^2}=\frac{n}{2},
 $$
 and the equality because the boundary has constant mean curvature $H=0$, besides the other restrictions cited in the results above.  Conversely, from \eqref{bD}, if $\psi$ is a (positive or negative) real Killing spinor on $\Omega$, we have
 $$
 {\bf D}\psi=\frac{n}{2}H\psi\pm\frac{n}{2}\gamma(N)\psi=\pm\frac{n}{2}\gamma(N)\psi.
 $$
 Now, we combine $\phi$ and $\gamma(N)\phi$ like in the statement of the theorem, taking into account \eqref{supercom} and we realise that $\frac{n}{2}$ it is, in fact, an  eigenvalue of ${\bf D}$ associated to the eigenspinor $\phi$ and, so, $-\frac{n}{2}$ associated to $\gamma(N)\phi$.
 \hfill $\square$  
 
\begin{remark}\label{Montiel-Ros}{\em 
Different to what happens with the Laplacian operator $\Delta^f$, acting on smooth functions, we can paraphrase the result above by asserting that, {\em with respect to the Dirac operator, all embedded minimal hypersurfaces in the sphere ${\mathbb S}^{n+1}$ are immersed by the first  $\pm\frac{n}{2}$-eigenvalues (see \cite{MR1}).} As for the Laplacian, we know that, for all the immersed minimal  hypersurfaces, the height functions are eigenfunctions corresponding to the eigenvalue $n$. Hence, $\lambda_1(\Delta^f)\le n$. But, in general, we do not know if  this eigenvalue $n$ is or not  the first eigenvalue of the operator $\Delta^f$.
{\em When the equality $\lambda_1(\Delta^f)= n$, we will show that the minimal surface $\Sigma$ is immersed in ${\mathbb S}^n$ by the first eigenvalues of the Laplacian} (see \cite{LiY}). It is known that the only metric on a $2$-dimensional sphere admitting a
minimal immersion into ${\mathbb S}^n$ by the first eigenfunctions is the standard one (this
follows, for example, from the fact that the multiplicity of the first eigenvalue
for such a metric is at most three, see the Cheng work \cite{Che}). In \cite{MR1}, we also  
showed that it is possible to extend this property for an arbitrary compact
surface, in the following way:
{\em 
For each conformal structure on a compact surface, there exists at most
one metric admitting a minimal immersion into a unit sphere by the first
eigenfunctions.}
As a consequence, the class consisting of such immersions seems not be too big.
This  enables us to characterise the equalities in some inequalities obtained by Li and Yau which relate the conformal area of a Riemannian
surface, the first non-zero eigenvalue of its Laplacian and the total mean
curvature for immersions of the surface  in the Euclidean sphere.
As the real projective plane has only one conformal structure, the only
metric on ${\mathbb R}{\mathbb P}^2$ admitting a minimal immersion into a sphere by the first
eigenfunctions is also the standard one. Thus, the metrics on ${\mathbb S}^2$ or  ${\mathbb R}{\mathbb P}^2$  having
this type of immersions are completely classified. Reasonably, we and some other authors have long been  interested
in extending this classification for other compact surfaces. Besides ${\mathbb S}^2$ and  ${\mathbb R}{\mathbb P}^2$,
the torus has the simplest family of conformal structures. The square and
flat equilateral  tori are the only known examples of Riemannian tori admitting
a minimal immersion into a sphere by the first eigenfunctions (these immersions lie in ${\mathbb S}^3$ and ${\mathbb S}^5$, respectively). For this surface we obtain the following
partial classification result \cite{MR1}:
{\em 
The only minimal torus immersed into ${\mathbb S}^3$  by the first eigenfunctions is the
Clifford torus. Moreover, we proved that 
there exist conformal structures on a torus for which there are no
metrics admitting a minimal immersion into any sphere by the first eigenfunctions.}

Note that this result puts before us two apparently encountered facts: it seems that it is considerably difficult to immerse minimally
a compact surface in a sphere by the first eigenvalues of its Laplacian, whereas we suspect that all the {\em embedded}
minimal surfaces in a sphere should have the same first eigenvalue $\lambda_1({\bf D})=\frac{n}{2}$. The
result quoted above gave us a relation between two well-known conjectures: the so-called Lawson's conjecture, which asserts
that the only torus minimally embedded into ${\mathbb S}^3$ is the Clifford torus, and the denominated
Yau's conjecture (or \#100 Problem in \cite{Y2}), which says that each minimal embedding of a compact surface
into ${\mathbb S}^3$ is by the first eigenfunctions of the Laplacian. From these results, we knew for a long time, that the Yau conjecture was
true, it would follow that two compact surfaces minimally
embedded into ${\mathbb S}^3$ are isometric provided that they are conformally equivalent and, so, the Lawson
conjecture should be also true.  Li and Yau estimate the conformal area of a torus
in terms of the area and the first eigenvalue of the only flat metric existing for
each conformal structure. We improve their bound for the conformal area and
this enables us to enlarge the family of conformal structures on a torus for
which the Willmore conjecture is satisfied. Since then (it was the beginning of the nineties 
of the past century), one of the two conjectures has recently become a nice theorem. Indeed, Brendle proved:
{\em
Any embedded minimal torus in ${\mathbb S}^3$ is congruent to the Clifford torus. This answers the conjecture posed by H.B. Lawson Jr. in 1970.
}
However, as far as we know, the famous {\em Problem \#100} in \cite[1982]{Y2} remains unresolved.  We fiercely encourage the reader to see \cite{Br3}, in particular Chapter 5. 
}
\end{remark}

 \begin{remark}{\em It is also worthy to note that, as another consequence of the estimate $\lambda_1({\bf D})\ge \min_\Sigma\sqrt{1+H^2}$ in Theorem \ref{lowerSS}, if $\Sigma$ is a compact boundary in the sphere ${\mathbb S}^{n+1}$  and we know that $\lambda_1({\bf D})\le \frac{n}{2}\sqrt{1+H^2}$, with $H\ge 0$, then we have the equality is  $H$ of $\Sigma$ is be constant on $\Sigma$. Moreover, $\Omega$ supports the existence of, at least, one non-trivial real Killing spinor (see list in [B\"a]). Hence, $\Omega$ has to be Einstein with positive scalar curvature. So, all our hypotheses are close $\Omega$ to be a spherical domain bounded by an embedded hypersurface with constant mean curvature. This would be the most similar to the solution to the original spherical Min-Oo conjecture and it is the third of our announced results as consequences of Theorem \cite{BMN,Min-Oo}.} 
\end{remark}

\begin{theorem}[Min-Oo's Conjecture]\label{FMOo1}
Suppose that $\Sigma$ is an enclosing hypersurface in an  $(n+1)$-dimensional compact connected manifold $\Omega$ with  scalar curvature $R\ge(n+1)$ with minimal boundary, that is $H\ge 0$. If $\Sigma$ admits no harmonic spinors and can be immersed isometrically and Riemannian isospinally (full) minimal immersion into a sphere ${\mathbb S}^{n+k}$, with $k\ge 1$, then $\Omega$ 
must be either a hemisphere  of  ${\mathbb S}^{n+1}\subset {\mathbb S}^{n+k}$  and $\Sigma$ its corresponding equator enclosing it. 
\end{theorem}
\pf\! 
The proof will follow the steps of that for Theorem  \ref{FMOo}. We combine Theorem \ref{lowerSS} and the spherical version of Remark \ref{highD} with \eqref{Hp} and, using that ${\overrightarrow H_0} =0$, we have 
\begin{equation}\label{IMM}
\frac{n^2}{4}\le \frac{n^2}{4}\min_\Sigma \left(1+{ H}^2\right) \le |\lambda_1(D)|^2\le \frac{n^2}{4}\frac{{\displaystyle\int_\Sigma \left(1+|{\overrightarrow {H_0}}|^2\right)}}{\hbox{\rm vol}(\Sigma)}= \frac{n^2}{4}.
\end{equation}
As we have obtained this chain of equalities, on the one hand, there a Killing spinor  defined on $\Omega$ and this compact domain lies in one the five types of spin manifolds listed in Theorem \ref{FMOo}. Moreover, according to the digression before the statement of that theorem, we have
$$
\frac{n}{2}H\gamma(N)\phi\mp\frac{n}{2}\gamma(N)\phi={\bf D} \phi=\frac{n}{2}\gamma(p+{\overrightarrow H})\phi=\frac{n}{2}\gamma(p)\phi.
$$
Taking norms, a first consequence is the relation $|{\overrightarrow H_0}|=H=0$, that is, $\Sigma$ is minimal in $\Omega$, as well. Moreover, $\gamma(N)=\pm\gamma(p)$, where $N$ is the inner unit normal in $\Omega$ and $p$ is the position vector in ${\mathbb S}^{n+k}\subset {\mathbb R}^{n+k+1}$. So, $k=0$ and ${\bar H}=HN\ge  0$ and we can apply Theorem {\bf B)}. Since the first eigenvalue $\lambda_1(D)$ consists of the restrictions to $\Sigma$ of the space of (positive or negative) real Killing spinors on $\Omega$, then its dimension is $N=(2^{\left[\frac{n}{2}\right]}, 2^{\left[\frac{n}{2}\right]})$, $(1,1)$, $(2,0)$, ($\frac{n+1}{4}+1$,0) or again $(1,1)$, according  to $\Omega$ lies in a sphere, in a Sasakian manifold, in a Einstein-Sasakian-manifold, in a 3-Sasakian manifold or, finally,  in a six-dimensional nearly-K\"ahler manifold.   On the other hand, this same eigenspace is $N=2^{\left[\frac{n+k}{2}\right]}$-dimensional, because we have said that it comes from the space of Killing spinor of the ambient ${\mathbb S}^{n+k}$. After a checking taking in account the table in \cite[p. 512]{Ba1}, all the cases can be reduced to two unique cases: 
either $(n+1,k)=(\hbox{\rm arbitrary},0)$ or $(n,k)=(\hbox{\rm arbitrary},1)$. In the first case, since the codimension $k=0$, $\Sigma$ is isometric to ${\mathbb S}^n$ and it is the minimal boundary of a domain $\Omega\subset {\mathbb S}^{n+1}$. As for the second case, the minimal boundary $\Sigma$ of the domain $\Omega\subset {\mathbb S}^{n+1}$ admits an isometric immersion in ${\mathbb S}^{n+1}$ but it does not admit harmonic spinors. But the  last one is impossible because of dimensional reasons (revisit the list in Theorem \ref{lambda1}). Thus, lastly, the unique remaining possibilities are that either $\Omega$ be a unit disc in ${\mathbb R}^{n+1}$ enclosed by $\Sigma={\mathbb S}^n$ or $\Omega$ is the domain enclosed in ${\mathbb S}^{n+1}\subset {\mathbb R}^{n+2}$ by any embedded minimal hypersurface $\Sigma$.  But the Gau{\ss} equation relating the scalar curvatures implies that $\Sigma$ is totally geodesic. Now, since $\Omega$ has a maximal number of Killing spinors and all their restrictions must provide a maximal number of Killing spinor on $\Sigma$, the unique valid solution is that of the hemisphere. \hfill $\square$  

\begin{remark}\label{upperH}{\rm
If we do not assume that the immersion $\psi: \Sigma^n\rightarrow {\mathbb S}^{n+k}$ is minimal, then we have from \eqref{IMM},
\begin{equation}\label{IMM1}
 0<\frac{n^2}{4}\min_\Sigma \left(1+{ H}^2\right) \le |\lambda_1(D)|^2\le \frac{n^2}{4}\frac{{\displaystyle\int_\Sigma \left(1+|{\overrightarrow {H_0}}|^2\right)}}{\hbox{\rm vol}(\Sigma)}.\nonumber
\end{equation}
This chain of inequalities is valid for any embedded hypersurface in ${\mathbb S}^{n+1}$ and any immersion  $\psi: \Sigma^n\rightarrow {\mathbb S}^{n+k}$ So, in these circumstances,
$$
{\hbox{\rm vol}(\Sigma)}\le \frac{{\displaystyle\int_\Sigma \left(1+|{\overrightarrow {H_0}}|^2\right)}}{\min_\Sigma \left(1+{ H}^2\right)},
$$ where $H$ and $\overrightarrow{H_0}$ are, respectively, the mean curvature function and the mean curvature vector of the initial embedding and the immersion $\psi$. Thus, the space of conformal immersions $\psi$ from $\Sigma$ to ${\mathbb S}^{n+k}$ whose mean curvature vector $|{\overrightarrow H}|$ and 
is bounded from above by a constant which only depends of  (the fixed) $  H$. Like in \cite{ChSch} by Choi and Schoen, we will fall into the temptation to pose a compactness conjecture: 
{\em 
Is it true that the set of conformal immersions of a  compact manifolds $\Sigma^n$, endowed with a fixed conformal structure, into ${\mathbb S}^{n+k}$ whose mean normal fields ${\overrightarrow H}$ have length bounded  by a fixed constant $C$ , $|{\overrightarrow H}|\le C$, is compact (in a suitable topology)?}
}
\end{remark}
  
\begin{corollary}[Well-posed (?) Original  (Spherical) Min-Oo's Conjecture]
Suppose that $\Sigma$ is an enclosing  hypersurface of an  $(n+1)$-dimensional compact spin Riemannian manifold  $\Omega$ with scalar curvature $R\ge n(n+1)$. If, moreover, $\Sigma$ is  isometric to a unit sphere ${\mathbb S}^{n+1}$ and is totally geodesic, then $\Omega$ must be a hemisphere ${\mathbb S}^{n+1}_+$ of a unit sphere ${\mathbb S}^{n+1}$ of radius $1$ and $\Sigma$ is the corresponding equator ${\mathbb S}^n$ enclosing it. 
\end{corollary}
\pf\!
 The sphere ${\mathbb S}^n$ has an only topologically spinorial structure which does not admit harmonic spinors (with whatever the Clifford product is) and that, evidently, it admits a totally 
 geodesic isometric immersion in ${\mathbb S}^{n+1}$ and, so  $H=0$. Hence, $\Omega$ is either a hemisphere of ${\mathbb S}^{n+1}$ and $\Sigma$ is the corresponding equator. 
 \hfill $\square$
\begin{remark}{\em 
Obviously the condition $R\ge n(n+1)$ is necessary, because otherwise one can 
perturb the hemisphere at an interior point so that $R\le n(n+1)-\varepsilon$, for some
small $\varepsilon > 0$, without changing the assumptions on the boundary \cite{HW1}. Moreover this  usual assumption $R\ge n(n+1)$ has been invalidated by the counterexamples built by Brendle, Marques and Neves \cite{BMN} and, very recently, by \cite{Sw}, since all those examples require at least one point with this strict inequality in the bulk manifold. {\em Hence, we have been brought to add some new hypothesis about the Riemannian spin structure of the boundary $\Sigma$}, namely, the {\em isospinallity} of the embedding of $\Sigma$ in $\Omega$. Note that ${\mathbb S}^6$ has two spinorial metric structures coming from the usual ${\mathbb S}^7$: the usual one and that coming from its nearly-K\"aler structure. They are different because the dimensions of their Killing spinor spaces. } 
\end{remark}

As in the flat case, we can obtain a kind of positivity for a possible quasi-local mass in this new context.

\begin{theorem}[Brown-York  in mean-convex and spherical case] 
Let $\Omega$ be a spin Riemannian manifold of dimension $n+1$ with  scalar curvature $R\ge n(n+1)$ and having a non-empty boundary $\Sigma$ whose inner mean curvature $H\ge 0$ is also non-negative (mean-convex) and without harmonic spinors. Suppose that there is an isometric and isospin immersion from $\Sigma$ into another spin manifold $\Omega_0$ carrying on a non-trivial real Killing spinor field and let $H_0$ its mean curvature with respect to any of its orientations.  Then, we have$$
\int_\Sigma\sqrt{1+ H^2}\le \int_\Sigma\sqrt{1+H_0^2},$$
provided that the boundary  does not admit harmonic spinors.
The equality implies that $H=H_0$. Then, if $n=2$, $\Omega_0$ is a domain in ${\mathbb S}^3$ and the two embeddings differ by a direct rigid motion.

\end{theorem}

{\bf Proof.}  Denote by $\psi$ the spinor on $\Omega_0$ and let $\phi=\psi_{|\Sigma}$ its restriction onto $\Sigma$ through the existent immersion. Let's recall that the parallelism of $\psi$ gives
$$
{\bf D}\phi=\frac{n}{2} \sqrt{1+H_0}\,\phi\qquad\hbox{and}\qquad |\phi|=1.
$$  
Now, we apply Theorem 2 and have the desired inequality$$
\int_\Sigma \sqrt{1+H^2}\le \int_\Sigma \sqrt{1+H^2_0}.$$

If the equality is attained, then 
$$
\frac{n}{2}\sqrt{1+H_0^2}={\bf D}\phi=\frac{n}{2}\sqrt{1+H^2}
$$ 
and so $H=H_0\ge 0$.
Then, the immersion of $\Sigma$ into the second ambient space $\Omega_0$ is mean-convex as well (with respect to
the choice of inner normal to $\Omega$).  

When $n=2$, from this equality and the fact that $K=K_\phi$ (because the two embeddings are isometric and preserve the Gauss curvatures), we deduce that the two second fundamental forms coincide. The Fundamental Theorem of the Local Theory of Surfaces allows us to conclude that the two boundaries differ by a direct rigid motion of the Euclidean space. 
 
\begin{corollary}[Cohn-Vossen Rigidity Theorem in the Sphere]
Two isometric and isospin mean-convex compact surfaces embedded in a sphere ${\mathbb S}^3$ without harmonic spinors must be congruent, provided they do not admit harmonic spinors.
\end{corollary}
{\bf Proof.} Let $\Omega$ and $\Omega_0$ be the two domains determined in ${\mathbb S}^3$ by two corresponding surfaces identified by means of an isometry. 
Then, we can apply Theorem 8 interchanging the roles of $\Omega$ and $\Omega_0$ and take in mind the case of the equality. \hfil\qed

The fact that the Cohn-Vossen is closely related with the total squared mean curvature of the submanifold already appeared in \cite[Theorem 7.18]{MR2}

\section{Conformal covariance an applications}
We have already made reference to the conformal covariance \eqref{confD} of the Dirac operator on any spin Riemannian manifold discovered by Hitchin in \cite[Section 1.4]{Ht},
as we recalled  in Remark \ref{confor} in Section \ref{NNC}.  Let $\Sigma$ be a spin Riemannian manifold and $\left<\,,\right>$ and ${\bar g}
=e^{2u}g=e^{2u}\left<\,,\right>$, where $u$ is a smooth function defined on the manifold,  two pointwise conformal metrics. It is well-known that one can identify the associated spinor bundles, the spinor metric and connection and the Clifford multiplication. But, the Dirac operators corresponding to the two metrics are different, though they are closely related. Indeed, given a spinor field $\psi\in\Gamma({\bf S}\Sigma)$, we define 
another one ${\bar \psi}$ by the relation
\begin{equation}\label{conffi}
{\bar \psi}=e^{-\frac{n-1}{2}u}\psi.
\end{equation}
Then, we have the following aforementioned conformal weighted covariance
\begin{equation}\label{confcov}
{\bar {\bf D}}({\bar \psi})={\bar{\bf D}}(e^{-\frac{n-1}{2}u})=e^{-\frac{n+1}{2}u}{\bf D}\psi,
\end{equation}
where, of course, ${\bar{\bf D}}$ and ${\bf D}$ are the Dirac operators corresponding to the metrics ${\bar g}$ and $\left<\,,\right>$, respectively. The reader may be interested in how to use this covariance to obtain some results relating $\lambda_1({\bf D})$ and conformal invariants such that the Yamabe invariant or the first eigenvalue of the Steklov problem on $\Sigma$. In this case, take a look at \cite{Hi,HMZ2}, for example. On the other hand, if $d\mu$ is the Riemannian measure corresponding to the metric $\left<\,,\right>$, it is clear that 
\begin{equation}\label{confmu}
d {\bar \mu}=e^{nu} d\mu,
\end{equation}
for the Riemannian measure of $g$. A straightforward computation from the three equalities above leads us to obtain 
\begin{eqnarray}\label{confiz}
&|{\bar{\bf D}}{\bar \psi}|\,|{\bar \psi}|\,d{\bar\mu}=e^{-\frac{n-1}{2}u}\,|{\bf D}\psi|\,e^{-\frac{n+1}{2}u}|\psi|\,e^{nu}\,d\mu=
|{{\bf D}}{\psi}|\,|{\psi}|\,d{\mu}&\\
&\sqrt{|{\bar{\bf D}}{\bar \psi}|^2|{\bar \psi}|^2-\left<{\bar{\bf D}}{\bar \psi},\gamma(N){\bar \psi}\right>^2}\,d{\bar\mu}=\sqrt{|{{\bf D}}{\psi}|^2|{ \psi}|^2-\left<{{\bf D}}{ \psi},\gamma(N){ \psi}\right>^2}\,d{\mu}\nonumber&.
\end{eqnarray}
This means that the left-hand side of our initial Theorems {\bf A} and {\bf B} is intrinsically conformally covariant under the change of the spinor $\psi$ for ${\bar\psi}$. As for the extrinsic right-hand side, let us suppose from now on that the function $u$ is defined on the whole $\Omega$ and the conformal change affects to the whole of the bulk. Thus, if we remind that the Willmore integrand $H^n\,d\mu$  
is conformally invariant, that is, that ${\bar H}^n\,d{\bar\mu}=H^n\,d\mu$. Then
\begin{equation}\label{confH}
{\bar H}=e^{-u}H
\end{equation}
and, as a consequence,
\begin{equation}\label{confder}
{\bar H}\,|{\bar \psi}|^2\,d{\bar\mu}=e^{-u}H\,e^{-(n-1)u}|\psi|^2\,e^{nu}\,d\mu={H}\,|{\psi}|^2\,d{\mu},
\end{equation}
that is, the right-hand side of  Theorem {\bf A} and Theorem {\bf B} is extrinsically conformally covariant under the same change of the spinor $\psi$ for ${\bar\psi}$
as above. 

Thus, let us suppose that {\em $\Omega$ is a compact domain of ${\mathbb S}^n$} and that we choose a metric  ${\bar g}=e^{2u}\left<\,,\right>$ conformal to the Riemannian structure $g=\left<\,,\right>$ induced on the boundary $\Sigma$ from the metric of the bulk manifold $\Omega$. As usual, let us suppose that the scalar curvature $R$ of $(\Omega,g)$ is greater than or equal to  $n(n+1)$, that the mean curvature $H$ is non-negative (see \eqref{confH}) and that $\Sigma$ has no harmonic spinors (note that, from \eqref{confcov}, this fact depends only of the conformal class $[g]$ spanned by $g$). Then, from the intrinsic conformal covariance of the left-hand side in \eqref{confiz} in Theorem {\bf B}, we deduce
$$
\int_\Sigma \sqrt{|{\bar{\bf D}}{\bar \phi}|^2|{\bar\phi}|^2-\left<{\bar {\bf D}}{\bar \phi},\gamma(N){\bar \phi}\right>^2}\,d{\bar\mu}\ge\frac{n}{2}\int_\Sigma H\,|\phi|^2\,d\mu
$$
for all spinor field $\phi:\Sigma\rightarrow {\mathbb C}^{{\left[\frac{n+1}{2}\right]}}$. Hence,
$$
\lambda_1({\bar{\bf D}})\int_\Sigma |{\bar\phi}|^2\,d{\bar\mu}\ge \frac{n}{2}\int_\Sigma H\,|\phi|^2\,d\mu.$$
But, using \eqref{conffi} and \eqref{confmu}, we have 
$$
|{\bar\phi}|^2\,d{\bar\mu}=e^{-(n-1)u}\,e^{nu}\,|\phi|^2\,d\mu.
$$
Thus,
$$
\lambda_1({\bar{\bf D}})\int_\Sigma |e^{{u}}\phi|^2\,d{\mu}\ge \frac{n}{2}\int_\Sigma \sqrt{1+H^2}\,|\phi|^2\,d\mu,
$$
for all spinor $\phi$ and all function $u$ defined on $\Omega$. Let us suppose that the function $u$ is non-positive, that is, $e^u\le 1$, pointwise on the boundary $\Sigma$. Geometrically, this means that the {\em new} metric ${\bar g}$ is less than or equal to the {\em old} one $g$, namely, ${\bar g}\le g$, pointwise as well. In other words, lengths of curves and distances between points are {\em contracted} by means of such a conformal change of metrics on $\Sigma$. Then, the inequality above becomes
$$
|\lambda_1({\bar{\bf D}})|\int_\Sigma |\phi|^2\,d{\mu}\ge \frac{n}{2}\int_\Sigma \sqrt{1+H^2}\,|\phi|^2\,d\mu,
$$
Hence, we have just proved the following result which is an improvement of Theorem \ref{lowerSS}. This improvement will lead us to the Llarull solution \cite[Theorems A, B and C]{Ll} (cf. \cite[Theorem 1.2]{LSW}) to a conjecture by Gromov (see reference was the already cited \cite{LSW}). 

\begin{theorem}\label{Llarull}
 Consider a compact spin  Riemannian manifold $\Omega$ of dimension $n+1$ whose metric $g$ has  scalar curvature $R\ge n(n+1)$  and its boundary $\Sigma$  is mean-convex without harmonic spinors. Then, if ${\bar g}$ is another metric on $\Sigma$ such that ${\bar g}\le g$,
 $$
\lambda_1({\bar {\bf D}})\ge\frac{n}{2}\min_\Sigma\sqrt{1+H^2}.
$$
The equality holds if and only if ${\bar g}=g$ and $(\Omega, g)$ supports a (positive or negative) real Killing $\psi$ which induces on $\Sigma$  a spinor field $\phi$ induced by $\psi$  (see \cite{Ba1,BFGK} and \cite[Appendix A]{Gi}),
by
$$
\phi=\psi\mp \left(H\pm\sqrt{1+H^2}\right)\gamma(N)\psi.
$$\hfill $\square$
\end{theorem}

As a almost immediate consequence, we obtain  the aforementioned Llarull result.

\begin{theorem}[Llarull's Solution to a Gromov's Conjecture, \cite{Ll,LSW}]\label{Llarull2}
Let  $g_0$ be the usual metric of the unit sphere ${\mathbb S}^n$ and ${\bar g_0}$ another metric on the sphere such that   ${\bar g_0}\le g_0$ pointwise. 
Then ${\bar g}_0=g_0$.
 \end{theorem}
\pf \, We put $\Omega={\mathbb S}^{n+1}_+$, a closed hemisphere of ${\mathbb S}^ {n+1}$ and consider on it the Euclidean metric $g_0$. Then, its boundary $\Sigma$ is the sphere ${\mathbb S}^n$ endowed with its canonical metric of constant curvature identically $1$. It is clear that $R=n(n+1)$ on $\Omega$ and that $\Sigma$ is mean-convex  (in fact, it is totally geodesic). Moreover, we know that ${\mathbb S}^n$ has  an only spinorial structure which does not admit harmonic spinors. Now, consider another metric ${\bar g_0}$ on the unit $n$-sphere such that ${\bar g_0}\le  g_0$ at each point. Since each two metrics on the sphere must be conformal up to a diffeomorfism, we can work by thinking that ${\bar g_0}=e^{2u}g$ for a certain function $u$ defined on ${\mathbb S}^n$ that we can extend to the whole of ${\mathbb B}^{n+1}$. The assumption ${\bar g_0}\le g_0$ is automatically translated into $u\le 0$.  Moreover, from \eqref{confH}, we know that ${\bar H}=0$, that is, $\Sigma$ is also minimal in $({\mathbb S}^{n+1}, {\bar g})$. We conclude taking into account Theorem \ref{Llarull} proved above and Remark \ref{upperH}.
\hfill $\square$

\begin{remark}{\rm 
Note that the reasoning in this Theorem \ref{Llarull2} remains to be valid when the ambient space ${\mathbb S}^{n}$  is replaced by a strictly convex hypersurface $\Sigma$ of ${\mathbb R}^{n+1}$. Indeed, in this case, $\Sigma$ is the boundary of an $(n+1)$-dimensional convex body $\Omega$ of ${\mathbb R}^{n+1}$ and we may endow it with the metric $g_0=N^*g_{{\mathbb S}^n}$, where $N:\Sigma \rightarrow {\mathbb S}^n$ is the (inner) Gau{\ss} aplication of the convex hypersurface which, in this case, is a diffeomorphism. This generalization of Llarull's result is due to Li, Su and Wang \cite[Theorem 1.2]{LSW} when $n$ is odd. 
}
\end{remark}

\end{document}